\magnification  1200\baselineskip13pt\vsize=8.7
truein\hsize6.5 truein\hoffset=-.2in 
\overfullrule=0pt

\def\Re{{ \rm {Re}\,}}\def\Im{  {\rm {Im}\,}}\def\sRe{\ss{\rm Re}\,}
\def\wtilde{\widetilde}

\def\E{{\cal E}}
\font\scaps=cmcsc10

\def\bar{\overline}
\font\smalsmalbf=cmbx8

\def\z{\zeta}


\font\smalltenrm=cmr8
\font\smallteni=cmmi8
\font\smalltensy=cmsy8
\font\smallsevrm=cmr6   \font\smallfivrm=cmr5
\font\smallsevi=cmmi6   \font\smallfivi=cmmi5
\font\smallsevsy=cmsy6  \font\smallfivsy=cmsy5
\font\smallsl=cmsl8      \font\smallit=cmti8

\def\smallfonts{\lineadj{80}\textfont0=\smalltenrm  \scriptfont0=\smallsevrm
                \scriptscriptfont0=\smallfivrm
    \textfont1=\smallteni  \scriptfont1=\smallsevi
                \scriptscriptfont0=\smallfivi
     \textfont2=\smalltensy  \scriptfont2=\smallsevsy
                \scriptscriptfont2=\smallfivsy
      \let\it\smallit\let\sl\smallsl\smalltenrm}


\font\smathbold=msbm5
\def\Ca{{\cal C}}\font\mathbold=msbm9 at 10pt\def\Hn{{\hbox{\mathbold\char72}}^n}\def\Sn{{S^{2n+1}}}
\font\smathbold=msbm7\def\sHn{{{\hbox{\smathbold\char72}}_{}^n}}
\def \C{{\hbox{\mathbold\char67}}}
\def\N{{\hbox{\mathbold\char78}}}\def\S{{\cal S} }
\def\R{{\hbox{\mathbold\char82}}}

\def\zbe{\z\!\cdot\bar\eta}\def\zbn{\z\!\cdot\bar\n}
\def\avg{-\hskip-1.1em\int}
\def\bze{\bar\zeta\!\cdot\!\eta}\def\zn{\z_{n+1}}\def\bz{{\bar {z}\,}}

\def\e{\epsilon}
\def\con{{\aut}}
\font\gothic=eusm8
\def\n{{\hbox{\gothic\char78}}}
\def\P{{\cal P}}\def\RP{{\R\P}}

\font\new=eurm10
\def\aut{{\hbox{\new\char65\new\char117\new\char116}}}


\def\lineadj#1{\normalbaselines\multiply\lineskip#1\divide\lineskip100
\multiply\baselineskip#1\divide\baselineskip100
\multiply\lineskiplimit#1\divide\lineskiplimit100}

\def\remark#1.{\medskip{\noin\bf Remark #1.\enspace}}
\def\endpf{$$\eqno/\!/\!/$$}

\def\pf#1.{\smallskip\noin{\bf  #1.\enspace}}

\def\noin{\noindent}

\def\ds{\displaystyle}
\def\ts{\textstyle}

\def\e{\epsilon}\def\part{\partial_t}

\def\isn{\int_{S^{2n+1}}}\def\wtilde{\widetilde}

\def\Rn{\R^n}\def\B{{\cal B}}
\def\sgn{{\rm sgn}}

\def\ref#1{{\bf{[#1]}}}

\def\p{\partial}
\def\omegan{\omega_{2n+1}}

\def\ints{\int_\Sigma}

\def\nhalf{{\ts{n\over2}}}\def\half{{\ts{1\over2}}}
\def\A{{\cal A}}\def\D{{\cal D}}\def\RR{{\cal R}}

\def\J{{\cal J}}
\def\T{{\cal T}}
\def \L{{\cal L}}

\def\ker{{\rm Ker}}

\def\GammaQ{\Gamma\big({Q\over2}\big)}

\def\H{{\cal H}}


\centerline{\bf{  Moser-Trudinger and Beckner-Onofri's inequalities on the CR sphere}}
\bigskip\centerline{Thomas P. Branson, Luigi Fontana, Carlo Morpurgo}
\vskip1.5em
\midinsert
{\smalsmalbf Abstract. }{\smallfonts 
We derive sharp Moser-Trudinger inequalities on the CR sphere. The first type is in the Adams form,
for powers of the sublaplacian and  for  general spectrally 
defined operators on the space of CR-pluriharmonic functions. We will then obtain the sharp 
 Beckner-Onofri inequality for CR-pluriharmonic functions on the sphere, and, as a consequence, a sharp logarithmic Hardy-Littlewood-Sobolev inequality  in the form given by Carlen and Loss.}

\endinsert
\vskip2em 
\centerline{\bf Table of contents}\bigskip
\item{0.} {\scaps   Introduction}\smallskip
\itemitem {-}{\sl Motivations and history}
\itemitem{-}{\sl Main results}
\itemitem{-} {\sl Ideas for related research}
\itemitem{-} {\sl In memory of Tom Branson}
\itemitem{-}{\sl Acknowledgments}
\itemitem{-}{\sl Addendum}\smallskip
\item{1.} {\scaps Intertwining operators on the CR sphere}\smallskip
\itemitem {-}{\sl The Heisenberg group, the complex sphere and the Cayley transform}
\itemitem {-} {\sl Sublaplacians on $\Hn$ and $\Sn$}
\itemitem {-}{\sl  Spherical and zonal harmonics on the CR sphere}
\itemitem{-} {\sl Hardy spaces and  CR-pluriharmonic functions}
\itemitem {-} {\sl Sobolev spaces}
\itemitem{-} {\sl Intertwining and Paneitz-type operators on the CR sphere}
\itemitem{-} {\sl Conditional intertwinors}
\itemitem{-} {\sl Intertwining operators on the Heisenberg group}
\itemitem{-} {\sl Intertwining operators and change of metric}
 
\smallskip
\item{2.} {\scaps Adams and Moser-Trudinger inequalities on the CR sphere}\smallskip
\itemitem{-} {\sl Adams inequalities for convolution type operators on the CR sphere}
\itemitem{-} {\sl Moser-Trudinger inequalities for operators of $d-$type on Hardy spaces}
\itemitem{-} {\sl Moser-Trudinger inequalities for powers of sublaplacians }

\smallskip
\item{3.} {\scaps Beckner-Onofri's inequality}
\smallskip
\item{4.} {\scaps The logarithmic Hardy-Littlewood-Sobolev inequalities}
\smallskip
\item{5.} {\scaps Appendix}
\bigskip
\centerline{\bf{0. Introduction}}\medskip \noin{\sl Motivations and history.}\medskip 
The problem of finding optimal Sobolev 
inequalities continues to be a source of inspiration to many analysts.
The literature on the subject is vast and rich.
Besides its intrinsic value, the determination of best constants in Sobolev, or 
Sobolev type, inequalities has almost always revealed or employed  deep facts 
about the geometric structure of the underlying space. More importantly, such 
constants  were often  the crucial elements needed   to
identify extremal geometries, and to solve important problems such as 
isoperimetric inequalities,  eigenvalue comparison theorems, curvature 
prescription equations, existence of solutions of PDE's, and  more. 

This kind of research has produced a wealth of conclusive results  in the context of 
Euclidean spaces and Riemannian manifolds.  In contrast, very 
little is known in subRiemannian geometry, even in the simplest cases of the 
Heisenberg group or the CR sphere; this is especially true with regards to best 
constants in Sobolev embeddings and sharp geometric inequalities.

In order to motivate our work, we present three by now classical 
sharp inequalities on the Euclidean   $\R^n$ and $S^n$. 
First, there is the standard Sobolev embedding $W^{d/2,2}\hookrightarrow 
L^{2n/(n-d)}$, $\,(0<d<n)$ represented by the   optimal inequality 
$$\|F\|_{q}^2\le C(d,n)\int_X F A_d F\qquad\quad q={2n\over n-d}\eqno (0.1)$$
with $C(n,d)=\omega_n^{-d/n}
\Gamma\big({n-d\over2}\big)/\Gamma\big({n+d\over2}\big)$, and where $\omega_n$ denotes the volume of $S^n$. For $X=\R^n$ the 
operator $A_d$ is $\Delta^{d/2}$, where $\Delta$ is the  positive Laplacian on $\R^n$,  and the extremals in (0.1) are dilations and 
translations of the function
$(1+|x|^2)^{-n/q}$. For $X=S^n$ the operator $A_d$ is the 
spherical picture  of $\Delta^{d/2}$, obtained from it via the stereographic 
projection and  conformal invariance.
These operators act on the $k$th order spherical harmonics 
$Y_k$ of $S^n$ as
$$A_d Y_k={\Gamma\big(k+{n+d\over2}\big)\over \Gamma\big(k+{n-d\over2}\big)} 
Y_k.\eqno(0.2)$$
When $d=2$, $A_2=\Delta_{S^n}+{n(n-2)\over4}$, the conformal Laplacian; 
for general  $d\in(0,n)$
$A_d$ is the  intertwining operator of order 
$d$ for the complementary series representations of $SO(n+1,1)$, and it is
an elliptic  pseudodifferential operator with the same leading symbol 
as $(\Delta_{S^n})^{d/2}$.
 The  fundamental solution of $A_d$ is given 
by the chordal distance 
 function $c_d|\zeta-\eta|^{d-n}$, with $\zeta,\eta\in S^n$, where $c_d$ is the 
same constant appearing in the fundamental solution (Riesz kernel) for 
$\Delta^{d/2}$ on $\R^n$. Higher order conformally  invariant powers of the Laplacian  on general manifolds were found  by Graham, Jenne, Mason, Sparling [GJMS], and are now known as the ``GJMS operators".
The  extremals for the inequality (0.1)  are the functions 
$|J_\tau|^{1/q}$ where $|J_\tau|$ denotes the density of the  volume change via a 
conformal transformation $\tau$ of $S^n$.
 
Inequality (0.1) is invariant under the action of the conformal group, 
both on $\R^n$ and $S^n$; for example on $\R^n$ in addition to the 
usual dilation/translation invariance,  there is also an invariance under 
inversion: the action $F\to F\big(x/|x|^2\big)|x|^{-2n/q}$ leaves both sides of 
(0.1) unchanged. It is this  aspect that makes this type of 
operators particularly   interesting.

For $d=2$ it was Talenti [Ta] who  first derived (0.1) on $\R^n$, followed by Aubin 
[Au1] on $S^n$. 
For general $d$ the inequality is the dual  of the sharp Hardy-Littlewood-
Sobolev inequality
obtained  by Lieb [L], a fundamental inequality which concerns the minimization 
of $\|F*|x|^{-\lambda}\|_q/\|F\|_p$ in the case $p=2$; as stated, (0.1) appears 
in [Bec].
\smallskip


Next, there is the limit case $d=n$ of (0.1), which  gives the so-called exponential class embedding
$W^{n/2,2}\hookrightarrow e^L$, and more generally $W^{d,n/d}\hookrightarrow e^L$, itself 
a limiting case of $W^{d,p}\hookrightarrow L^{np/(n-dp)}$. 
In concrete terms the Sobolev embedding in the critical case $dp=n$  is represented by an Adams-Moser-Trudinger inequality of type

$$\int_{X} \exp\bigg[\alpha_{d}\bigg({|F|\over \|B_{d}F\|_p}\bigg)^{p'}\bigg]\le 
c_0,\qquad p={n\over d}\eqno(0.3)$$

\noin where  $B_{d}$ is a suitable, possibly vector-valued,  pseudodifferential operator of order $d$, and 
where
the constant $\alpha_{d}$ is best, i.e. it cannot be replaced by a larger
constant. Here $F$ runs through an appropriate subspace of $W^{d,n/d}$ where $B_{d}F\neq0$.

In the case of bounded domains of $\Rn$ the first sharp result is due to Moser [Mo1], who obtained (0.3) with best constant in the case $d=1$ and $B_1=\nabla$, for $F\in W_0^{1,n}(\Omega)$. Earlier, Trudinger [Tr]
proved a similar inequality, without best constant, with $\|\nabla F\|_n$  replaced by $\|F\|_n+\|\nabla F\|_n$, for $F\in W^{1,n}(\Omega)$.
 Adams [Ad] found the sharp version of Moser's inequality for   higher order gradients $B_k=\nabla^k$, and $F\in W_0^{k,n/k}(\Omega)$.  A few years later,  Fontana [F] extended Moser's and Adams' results to arbitrary compact manifolds without boundary. Since then, and up to recent times,  many authors worked on other extensions and generalizations of Moser's result, often motivated by problems in conformal geometry and nonlinear PDE.


\smallskip
On the sphere there is also another  sharp form of the exponential class embedding for $W^{n/2,2}(S^n)$, namely 
the so-called Beckner-Onofri inequality
$${1\over 2n!}\avg_{S^n} \!uA_nu+\avg_{S^n} \!u-\log\avg_{S^n} \!e^{u}\ge 0\ 
,\eqno(0.4)
$$
where $\ds{\avg}$ denotes the average operator, and where $A_n$ is the intertwining  operator 
defined by (0.2) in the limit case $d=n$ , with eigenvalues 
$\,k(k+1)...(k+n-1)$. Such $A_n$ is sometimes referred to as the Paneitz 
operator on the sphere, in honor of S. Paneitz who first discovered a fourth order conformally invariant operator on general manifolds, which reduces to $A_4=(\Delta_{S^4})^2+2\Delta_{S^4}$ on the Euclidean $S^4$. Note that  $A_2=\Delta_{S^2}$, when $n=2$. Due to the particular nature of $A_n$, the functional in (0.4) is invariant under the  group action $F\to F\circ\tau+\log|J_\tau|$, 
where $\tau$ is a conformal transformation of $S^n$; this action preserves the exponential integral. 
This important inequality was first derived by Onofri  in dimension 2, and its general $n-$dimensional form was discovered later by 
Beckner [Bec], via an endpoint differentiation argument based on (0.1) and  the sharp Hardy-Littlewood-Sobolev inequality. Later, Chang and Yang [CY] gave an alternative proof of (0.4) by a completely different method, based on an extended and refined version of the original compactness argument used by Onofri.

Estimate (0.4) has relevant applications in spectral geometry and mathematical physics, from comparison theorems for functional determinants to the theory of   isospectral surfaces. [Br], [BCY], [CY], [CQ], [O], [OPS]. 

Over the past couple of decades there has been a growing interest 
in finding 
the analogues  of the above results in the context of CR geometry. The biggest 
motivations are certainly the isoperimetric inequality, the isospectral problem, 
extremals for spectral invariants such as the functional determinant, and 
several other eigenvalue comparison theorems.

In the  CR setting, the first and only known sharp Sobolev embedding estimate of type (0.1) with conformal invariance properties
 is due to 
Jerison and Lee [JL1],[JL2], and it holds on the Heisenberg group $\Hn$ and on the CR sphere $S^{2n+1}$ in the 
case  $d=2$,  for the CR invariant Laplacian (which is the standard sublaplacian in the case of $\Hn$). 
 The corresponding version for operators of order  $0<d<Q=2n+2,\,d\neq 2$,  is only conjectured $\!\!\!{}^\dagger$\vfootnote\dag{
\hglue-.5em See the ``Addendum" at the end of this section, about a recent breakthrough made by Frank and Lieb [FL] in this regard.}, and involves the 
intertwining operators $\A_d$  for the complementary series representations of $SU(n+1,1)$. The explicit form of such operators has been known for quite some time, for example by  work of Johnson and  Wallach [JW], and also Branson, \'Olafsson and \O rsted [BO\O], and can be described as follows. Let $\H_{jk}$ be the space of harmonic polynomials of bidegree $(j,k)$
 on $\Sn,$ for  $j,k=0,1,...$; such spaces make up for the standard decomposition of $L^2$ into $U(n+1)$-invariant and irreducible subspaces. The intertwining operators of order $d<Q$ are characterized (up to a constant) by their action on $Y_{jk}\in\H_{jk}:$
$$\A_d Y_{jk}=\lambda_j(d)\lambda_k(d)Y_{jk},\qquad \lambda_j(d)={\Gamma\big(j+{Q+d\over4}\big)\over\Gamma\big(j+{Q-d\over4}\big)};\eqno(0.5)$$
when $d=2$ this gives the CR invariant sublaplacian.
As it turns out these operators have a simple fundamental solution of type $c_d|1-\zbe|^{d-Q\over2}$, where $\z,\eta\in \Sn$, for a suitable constant $c_d$. The conformally invariant  sharp Sobolev inequality that is conjectured to be true is 
$$\bigg(\avg_{\Sn} |F|^q\bigg)^{2/q}\le {1\over\lambda_0(d)^2}\,\avg_\Sn F \A_d F\qquad\quad q={2Q\over Q-d}\eqno (0.6)$$
with extremals  $|J_\tau|^{1/p}$, $\tau$ a conformal transformation of $\Sn$; this is the Jerison-Lee inequality for $d=2$ but it is an open problem for general $d$. This conjecture does not seem to appear on any published articles, but it is well known within the group of researchers interested in this type of questions. One of the aspects that makes the CR treatment more difficult,  is the lack, to date, of an effective symmetrization technique on the CR sphere or the Heisenberg group, that would allow for example to show the dual version of (0.6), namely the CR Hardy-Littlewood-Sobolev inequality.

Regarding Moser-Trudinger inequalities at the borderline case $d=Q/p$,  Cohn and Lu recently made some progress 
[CoLu1], [CoLu2],  deriving the CR analogue of  (0.3) with sharp exponential constant in the case of the gradient, $p=Q$, both on $\Hn$ or the 
CR $S^{2n+1}$ (see also [BMT] for similar results on Carnot groups). In regards to the ``correct" CR analogue of Beckner-Onofri's inequality  (0.4), the situation is not so obvious. One would certainly start to consider the operator $\A_Q=\lim_{d\to Q}\A_d$, the intertwining or Paneitz operator at the end of the complementary series range; the kernel of this operator is the space of CR-pluriharmonic functions on $\Sn$, given by $\P:=\bigoplus_{j>0}(\H_{j0}\oplus\H_{0j})\oplus\H_{00}$. On the basis of (0.4) the natural conjecture would be that  for a suitable constant $c_n$
$$c_n\avg F\A_QF  -\log\avg e^{F-\pi F}\ge0,\,\qquad \forall F\in W^{Q/2,2}\eqno(0.7)$$
where $\pi F$ denotes the Cauchy-Szego projection of $F$ on the space $\P$. The Euclidean version (0.4) can be cast in a similar form, with $\pi F$ being just the average of $F$. This inequality however  is {\sl not} invariant under the conformal action
that preserves the exponential integral, i.e. $F\to F\circ\tau+\log|J_\tau|$. On the other hand, the fact that   $\A_Q$ has such large kernel $\P$
combined with the invariance of $\P$  under the conformal action (see Prop. 3.2) leads one to think that there should be a CR version of $(0.4)$ that  is conformally invariant and whose natural  ``milieu" is the space of CR-pluriharmonic functions; in this work we show that this is indeed the case.

\bigskip\noin{\sl Main results.}\medskip

The CR version of Beckner-Onofri's inequality proven in this paper is described as follows.
Let $\A'_Q$ be the operator acting on CR-pluriharmonic functions as 
$$\A_Q' \sum_j(Y_{j0}+Y_{0j})=\sum_j\lambda_j(Q)(Y_{j0}+Y_{0j}),\qquad \lambda_j(Q)=j(j+1)..(j+n)$$
where $Y_{j0}\in\H_{j0},\,Y_{0j}\in\H_{0j}$. In Theorem 3.1 we prove that  for any real  $F\in W^{Q/2,2}\cap \P$ we have
$${1\over2(n+1)!}\avg_{\Sn}F\A_Q'F+\avg _{\Sn}F-\log\avg_{\Sn} e^F\ge0.\eqno(0.8) $$
The functional in  (0.8)  is invariant under the conformal action $F\to F\circ\tau+\log|J_\tau|$, where $\tau$ is a conformal transformation of $S^{2n+1}$ (i.e.  $\tau$ is identified with an element of $SU(n+1,1)$), and $|J_\tau|$ its Jacobian determinant. The extremals of (0.8) are precisely the  functions $\log|J_\tau|$.

A few remarks are in order. First, the conformal  action is  an {\it affine} representation of $SU(n+1,1)$, and the minimal nontrivial closed (real) subspace of $L^2$ that is invariant under such action is precisely the space of real CR-pluriharmonic functions (Prop. 3.2).
This is in contrast with  the Euclidean case, for the action induced by $SO(n+1,1)$, since  in that case
the only invariant closed subspaces of $L^2$ are the trivial ones. This observation seems to justify (at least partially) that inequality (0.8) could be regarded as the direct  CR analogue of (0.4) from the point of view of conformal invariance.

Secondly, the key character in (0.8) is the operator $\A_Q'$, which we call the {\it conditional intertwinor} of order $Q$ on $\P$. 
 This operator is the CR analogue on $\P$ of the Paneitz, or GJMS,  operator $A_n$ on the standard Euclidean sphere, and coincides, up to a multiplicative constant, with the {\it $d-$derivative} at $d=Q$ of $\A_d$ restricted to $\P$. Moreover, we have 
$$\A_Q'F=\prod_{\ell=0}^n\big({\ts{2\over n}}\L+\ell\big)F\,,\qquad F\in\P$$
where $\L$ is the standard sublaplacian on the sphere. To our knowledge such operator is introduced  here for the first time.

Finally, if conjecture (0.6) were true then (0.8) would result  by the same endpoint differentiation argument used by Beckner to obtain (0.4). The meaning of this is that even though  we do not know whether (0.6) holds, we can still consider the functional 
$$\J_d[G]={1\over\lambda_0(d)^2}\avg G\A_d G-\bigg(\avg|G|^{q}\bigg)^{2/q},\qquad\quad q={2Q\over Q-d}$$
and take the $d-$derivative at $Q$ of $\J_d[1+(1/q)F]$ under the restriction $F\in\P$; the result of this operation is the functional in (0.8).
This argument will in fact be used to prove the conformal invariance of (0.8) (see Prop. 3.2).

Our proof of (0.8) follows the same general strategy used by Chang-Yang and Onofri.
The first step is to show that the functional in (0.8) is bounded below. This is accomplished by a ``linearization" procedure from a sharp Adams/Moser-Trudinger inequality on the CR sphere derived here for the first time. Indeed, a  portion of this work is dedicated to inequalities of type 
$$\int_{\Sn}\exp\bigg[ A_{d}\bigg({|F|\over\|B_dF\|_p}\bigg)^{p'}\bigg]d\z\le C_0\eqno(0.9) $$
where $0<d<Q$, $dp=Q$, which are of independent interest. We will obtain (0.9) for what we call $d-$type operators on Hardy spaces $\H^p$, or $\P^p$ ($L^p$ boundary values of pluriharmonic functions on the ball), and which are  essentially finite sums of powers of the sublaplacian, restricted to such spaces, with  leading power equal to  $d/2$. When  $p=2$, the case of interest for (0.8), we have 
$A_{Q/2}={1\over2} (n+1)!\omega_{2n+1}$
 and this constant  is sharp, i.e. in (0.9) it cannot be replaced by a larger constant. We will also obtain (0.9) on the full $W^{d,Q/d}$ for $B_d=\L^{d/2}$ or $B_d=\D^{d/2}$, where $\L$ is the sublaplacian of the CR sphere, and $\D=\L+{n^2\over4}$ is the conformal sublaplacian, with sharp constants for any $d<Q$.
All of these  inequalities will be applications of recent results by Fontana and Morpurgo [FM], on Adams inequalities in a measure-theoretic setting; their   proofs will follow  from  asymptotically sharp growth estimates on the fundamental solutions of the operators $B_d$, in terms of their distribution functions.

The second step toward a proof of (0.8) is to establish  that the functional has a minimum, via a compactness argument based on an Aubin's type inequality. This inequality is essentially saying that if a function $F$  has vanishing center of mass, then an inequality like (0.8) holds with and  improved constant on the leading order term, but with added lower order terms.

The final step is a version of the argument given by Chang-Yang [CY] based on Hersch's old results [H], in order to characterize the extremals. As a byproduct we will obtain sharp inequalities for the first eigenvalue of $\A_Q'$  under conformal change of contact structure on $S^{2n+1}$.

In the final part of the paper we will derive  from (0.8) the following sharp  logarithmic Hardy-Littlewood-Sobolev 
inequality:
$$(n+1)\avg\avg\log{1\over |1-\zbe|}\,G(\z) G(\eta)d\z d\eta\le \avg G\log G\,d\z\eqno(0.10)$$
valid for all $G\ge0$ with the right hand side finite, and  ${\ds{\avg }}G=1.$
The inequality is conformally invariant under the action $G\to (G\circ\tau)|J_\tau|$, and its extremals are the functions $|J_\tau|$, with $\tau$ any conformal transformation. The logarithmic kernel in (0.10) is a fundamental solution of $\A'_Q$ as an operator acting on CR-pluriharmonic functions with  mean 0:
$$(\A_Q')^{-1}(\z,\eta)=-{2\over\GammaQ\omega_{2n+1}}\,\log{ |1-\zbe|.}$$

In the Euclidean context (0.10) was obtained by Carlen and Loss [CL] from the sharp inequality (0.1), cast in its dual form,  via endpoint differentiation.
In some precise sense (0.10) and (0.8) are dual of one another. Finally, we will derive an equivalent version of (0.10) 
on the Heisenberg group, using the conformal invariance of such inequality.

\bigskip

\noin {\sl Ideas for related research}
\medskip
The inequality obtained by Beckner and Onofri turned out  to be central in 
the problem of finding extremal geometries for the functional determinant of certain operators  on  compact Riemannian manifolds.
We expect the same to be true in the case of CR geometry, namely that an explicit computation of functional determinants
of conformally invariant operators, at least in low dimensions, would involve the functional in (0.8), and that (0.8) itself would be useful in solving extremal problems. 

At the dual end, the third author has shown  in [M1]  that the logarithmic Hardy-Littlewood-Sobolev inequality on $S^n$ was 
the analytic expression of an extremal problem for the regularized zeta function of the Paneitz operators. 
Likewise we expect the same to be true on the CR sphere. 

We hope that the results presented in this paper will serve as an  incentive to pursue these  matters, and in particular to motivate the explicit calculation of  functional determinants for low dimensional CR-manifolds.

\bigskip
\noin{\sl In memory of Tom Branson.}\medskip

Tom Branson wrote: {\sl ``What I have in mind is to generalize Beckner's sharp, invariant\break 
Moser-Trudinger inequality on $S^n$, which is a fact about conformal
geometry, to a fact about CR geometry, and eventually other rank 1
and higher rank geometries"} [Br1]. Chang and Yang gave an alternative, symmetrization-free proof of Beckner's inequality on  $S^n$; \break it was Branson's idea that we might attempt to {\sl``play the same game"} on the CR sphere. {\sl ``This is not just any example; it's the one people will be by far most interested in, because\break  of CR geometry"} [Br1].  The present  paper is the result of our efforts to prove that Tom Branson's original intuition  was indeed correct: yes, we can play the same game, but on the space of CR-pluriharmonic functions (and with considerably more difficulties). 

Tom Branson suddenly passed away in March 2006.
\bigskip

\noin{\sl Acknowledgments.}
\medskip

The authors would like to thank  Francesca Astengo,  Bill Beckner, Arrigo Cellina, Bent \O rsted, Marco Peloso, Fulvio Ricci and Richard Rochberg for helpful comments.

\bigskip
\noin{\sl Addendum.} \medskip After this work was completed an important and remarkable breakthrough was made by R. Frank and E. Lieb  [FL], who were able to prove the sharp Hardy-Littlewood-Sobolev inequality on $\Hn$, or its equivalent version (0.6) on $S^{2n+1}$. Their proof is symmetrization-free. The proof of the existence of the optimizers is based on a  sophisticated compactness argument, whereas the characterization of the extremals  is accomplished by a  clever  enhanced version of a Hersch type argument used originally  by Chang-Yang in [CY], and adapted to the CR setting in the present paper (see section 3). 

\eject
\bigskip\centerline{\bf {1. Intertwining operators on the CR sphere.}}
\bigskip
\noin{\sl The Heisenberg group, the complex sphere and the Cayley transform}

\medskip
The Heisenberg group $\Hn$ is $\C^n\times \R$ with elements $u=(z,t)$, $z=(z_1,...,z_n)$, and with 
group law
$$(z,t)(z',t')=(z+z',t+t'+2\Im z\cdot \bar z')$$
where we set $z\cdot \bar w=\sum_1^n z_j\bar w_j$, for $w=(w_1,...,w_n)$. The Lebesegue-Haar measure on $\Hn$ is denoted by $du$.

Throughout the paper we will often use the standard notation for the homogeneous dimension of $\Hn$:
$$Q=2n+2.$$

The sphere $S^{2n+1}$ is the boundary of the unit ball $B$ of $\C^{n+1}$. In  coordinates, $\z=(\z_1,...,\z_{n+1})\in S^{2n+1}$ if and only if   $\z\cdot\bar\zeta=\sum_1^{n+1}|\z_j|^2=1$. The standard Euclidean volume element of $S^{2n+1}$ will be denoted by $d\zeta$. 

\smallskip The  Heisenberg group and the sphere are equivalent via the Cayley transform $\Ca:\Hn\to S^{2n+1}\setminus(0,0,...,0,-1)$  given by 
$$\Ca(z,t)=\Big({2z\over 1+|z|^2+i t},{1-|z|^2-i t\over 1+|z|^2+i t}\Big)$$
and with inverse
$$\Ca^{-1}(\z)=\Big({\z_1\over1+\z_{n+1}},...,{\z_n\over1+\z_{n+1}},Im {1-\z_{n+1}\over
  1+\z_{n+1}}\Big).$$ 
We will use the notation
$$\n=\Ca(0,0)=(0,0,...,1).$$
The Jacobian determinant (really a volume density) of this transformation is given by 
$$|J_\Ca(z,t)|={2^{2n+1}\over \big((1+|z|^2)^2+t^2\big)^{n+1}}$$
so that 
$$\int_{S^{2n+1}} Fd\zeta=\int_{\sHn} (F\circ\Ca)|J_\Ca| du$$

The homogeneous norm on $\Hn$ is defined by 
$$|(z,t)|=(|z|^4+t^2)^{1/4}$$
and the distance from   $u=(z,t)$ and  $v=(z',t')$ is 
$$d((z,t),(z',t')):=|v^{-1}u|=\big(|z-z'|^4+(t-t'-2\Im(z\cdot \bar
z'))^2\big)^{1/4}$$

On the sphere the distance function is defined as 
$$d(\z,\eta)^2:=2|1-\z\cdot \bar\eta|=\big|\,|\z-\eta|^2-2i
\,\Im(\z\cdot\bar\eta)\big|=
\big(|\z-\eta|^4+4\cdot{\rm \Im}^2(\z\cdot\bar\eta)\big)^{1/2}$$
and a simple calculation shows that if $u=(z,t),v=(z',t')$, and $\z=\Ca(u),\,\eta=\Ca(v)$.
then
$${|1-\z\cdot \bar\eta|\over
  2}=|v^{-1}u|^2\big((1+|z|^2)^2+t^2\big)^{-1/2}\big((1+|z'|^2)^2+(t')^2\big)^{-1/2}\eqno(1.1)$$
i.e.
$$d(\z,\eta)=d(u,v)\bigg({4\over (1+|z|^2)^2+t^2}\bigg)^{1/4}\bigg({4\over (1+|z'|^2)^2+(t')^2}\bigg)^{1/4}.\eqno(1.2)$$
\bigskip
\noin{\sl Sublaplacians on $\Hn$ and $S^{2n+1}$.}
\medskip
The sublaplacian on $\Hn$ is the second order differential operator
$$\L_0=-{1\over 4}\sum_{j=1}^n (X_j^2+Y_j^2) $$
where $X_j=\ds{{\p\over\p x_j}+2y_j{\p\over\p t},\; Y_j={\p\over\p y_j}-2x_j{\p\over\p t}},$ and $\ds{\p\over\p t}$ denote the basis of the space of left-invariant vector fields on $\Hn$.  One can check that 
$$\L_0=-{1\over2}\sum_{j=1}^n(Z_j\bar Z_j+\bar Z_jZ_j) $$
where $$Z_j=\ds{{\p\over \p z_j}+i\bar z_j{\p\over \p t},\quad\bar Z_j={\p\over \p\bar z_j}-i z_j{\p\over \p t}}$$ and with $\ds{{\p\over \p z_j}={1\over2}\Big({\p\over\p x_j}-i{\p\over\p y_j}\Big),\;{\p\over \p\bar z_j}={1\over2}\Big({\p\over\p x_j}+i{\p\over\p y_j}}\Big)$.
\medskip
The fundamental solution of $\L_0$ was computed by Folland [Fo1]
and 
$$\L_0^{-1}(u,v)=C_2\,d(u,v)^{2-Q},\qquad C_2={{2^{n-2} \Gamma({n\over2})^2}\over\pi^{n+1}}$$
so that 
$$G(u)=\int_{\sHn}C_2|v|^{2-Q}F(v^{-1}u)dv=\int_{\sHn}\L_0^{-1}(u,v)F(v)dv$$
solves $\L_0 G=F$.
\medskip
On the standard sphere, the sublaplacian is defined similarly  as 
$$\L=-{1\over2}\sum_{j=1}^{n+1}(T_j\bar T_j+\bar T_jT_j) $$
where 
$$T_j={\p\over\p\zeta_j}-\bar\zeta_j \RR,\qquad\RR=\sum_{k=1}^{n+1} \z_k{\p\over\p \z_k},\eqno(1.3)$$
and where the $T_j$  generate the holomorphic tangent space $T_{1,0}S^{2n+1}=T_{1,0}\C^{n+1}\cap \C TS^{2n+1}$. Explicitly 
$$\L=\Delta+\sum_{j,k=1}^{n+1}\z_j\bar\zeta_k{\p^2\over\p \z_j\p \bar \z_k}+{n\over2}(\RR+\bar \RR)\eqno(1.4)$$
with $\Delta=-\sum_{j}\ds{\p^2\over\p \z_j \p \bar \z_j}$.  The trasversal direction is the real vector field
$$\T={i\over2}(\RR-\bar \RR)={i\over2}\sum_{j=1}^{n+1} \bigg(\z_j {\p\over \p \z_j}-\bar\zeta_j {\p\over \p \bar \z_j}\bigg)\eqno(1.5)$$
and $\C TS^{2n+1}$ is generated by the $T_j,\bar T_j,\T$.
\medskip
The {\sl conformal sublaplacian} on the sphere is defined as 
$$\D=\L+{n^2\over4}.$$
The fundamental solution of  $\D$ has been computed by Geller [Ge] (Thm. 2.1 with $\alpha=0$ and modulo volume normalization)
$$\D^{-1}(\z,\eta)=c_2 \,d(\z,\eta)^{2-Q},\qquad c_2={{2^{n-1} \Gamma({n\over2})^2}\over\pi^{n+1}}=2C_2\eqno(1.6)$$
in the sense that for smooth $F:S^{2n+1}\to\C$ the function 
$$G(\z)=\D^{-1}F(\z)=\isn c_2d(\z,\eta)^{2-Q}F(\eta)d\eta$$
satisfies $\D G=F$.

The peculiarity of $\D$ is its  direct relation with  $\L_0$ via the Cayley transform:
$$\L_0\Big((2|J_\Ca|)^{Q-2\over2Q}(F\circ\Ca)\Big)=(2|J_\Ca|)^{Q+2\over2Q}(\D F)\circ\Ca\eqno(1.7)$$
which can be readily established by using the explicit formulas for the fundamental solutions and (1.2). The multiplicative factor 2 in the above formula appears because we use the standard volume elements for $\Hn$ and $S^{2n+1}$ instead of the volume forms associated with the standard contact forms $\theta_0$, and $\theta$ of these two spaces. In this case indeed we have that 
$$\int_{\sHn} f\,\theta_0\wedge d\theta_0...\wedge d\theta_0=2^{2n}n! \int_{\sHn} f du=\isn F\,\theta\wedge d\theta...\wedge d\theta=2^{2n+1}n! \isn Fd\zeta$$ where $f=(F\circ\Ca)(2|J_\Ca|)$ (see Jerison-Lee [JL1)]. This also accounts for  the factor  2 in the relation $c_2=2C_2$. 
\bigskip

\noin{\sl Spherical and zonal harmonics on the CR sphere}.
\bigskip
 The  space $L^2(S^{2n+1})$, endowed with the inner product
$$(F,G)=\isn F\bar G\,d\zeta,$$
can be decomposed as $L^2(S^{2n+1})=\ds{\bigoplus_{j,k\ge0} \H_{jk}}$, where $\H_{jk}$ is the space of harmonic polynomials on $\C^{n+1}$ that are homogeneous of degree $j,k$ in the $\z$'s and $\bar\zeta$'s respectively, and restricted to the sphere.  The dimension of $\H_{jk}$ is 

$$\dim(\H_{jk})=m_{jk}:= {(j+n-1)!(k+n-1)!(j+k+n)\over n!(n-1)!j!k!}\eqno(1.8)$$
and if $\{Y_{jk}^\ell\}$ is an orthonormal basis of $\H_{jk}$ then
the zonal harmonics are defined as  
$$\Phi_{jk}(\z,\eta)=\sum_{\ell=1}^{m_{jk}} Y_{jk}^\ell(\z)\bar {Y_{jk}^\ell(\eta)}$$
The $\Phi_{jk}$  are  invariant under the transitive action of $U(n)$ and it turns out that
$$\Phi_{jk}(\z,\eta)=\Phi_{jk}(\zbe):={(j+n-1)!(j+k+n)\over \omega_{2n+1}n!j!}(\zbe)^{j-k}P_k^{(n-1,j-k)}(2|\zbe|^2-1)\eqno(1.9)$$
if $k\le j$, and 
$\Phi_{jk}(\z,\eta)=\Phi_{jk}(\zbe):=\bar{\Phi_{kj}(\zbe)}$, if $j\le k$, where $P_k^{(n,\ell)}$ are the Jacobi polynomials (see [VK], Section 11.3.2).
\medskip
In particular, since $P_0^{(n-1,j)}\equiv 1$ we have also

$$\Phi_{j0}(\zbe)={(j+n)!\over j!n! \omega_{2n+1}}\,(\zbe)^j={\Gamma\big(j+{Q\over2}\big)\over\Gamma(j+1)\Gamma\big({Q\over2}\big)\omega_{2n+1}}\,(\zbe)^j\eqno(1.10)  $$
 and $\Phi_{0k}(\zbe)=\bar{\Phi_{k0}(\zbe)}=\Phi_{k0}(\bze) $.

If $F\in L^2$ then 
$$F(\z)=\sum_{j,k\ge0}\isn F(\eta)\Phi_{jk}(\zbe)d\eta$$
the series being convergent in $L^2$.

\bigskip
\noin{\sl Hardy spaces and  CR-pluriharmonic functions.}
\medskip
\noin In the sequel we will use the following notations\medskip
$\H=\bigoplus_{j\ge0}\H_{j0}=\big\{ L^2 {\hbox { boundary values of holomorphic functions on the unit ball }}\big\}$\medskip
$\bar \H=\bigoplus_{j\ge0}\H_{0j}=\big\{ L^2 {\hbox { boundary values of antiholomorphic functions on the unit ball }}\big\}$\medskip
$\P=\bigoplus_{j>0}(\H_{j0}\oplus\H_{0j})\oplus\H_{00}=\big\{ L^2 {\hbox { CR-pluriharmonic functions }}\big\}$\medskip
$\RP=\big\{ L^2 {\hbox { real-valued CR-pluriharmonic functions }}\big\}$\medskip
$\H_0,\,\bar \H_0,\,\P_0,\,\RP_0=$ functions in $\H,\,\bar\H,\,\P,\,\RP$ with 0 mean. Note that $\H_{00}$ is the space of constant functions.\medskip 

The space $\H$ is the classical Hardy space for the boundary of the unit ball of $\C^{n+1}$.
The Cauchy-Szego projection from $L^2(\Sn)$ to $\H$ is given by the Cauchy-Szego kernel 
$$K(\z,\eta)={1\over \omega_{2n+1}(1-\zbe)^{n+1}}=\sum_{j\ge0}\Phi_{j0}(\zbe).$$
 The projection operator on $\P$
$$\pi:L^2(\Sn)\to\P$$ has kernel $2\Re K(\z,\eta)-\ds{1\over\omega_{2n+1}}$.
Denote by $\P^\perp$ the orthogonal complement of $\P$, with respect to the standard Hermitian product $\zbe$,  i.e.
$$L^2(\Sn)=\P\oplus \P^\perp.$$

\bigskip
The Hardy spaces for $p>1$ are defined similarly. $\H^p$ will denote the $L^p$ closure of boundary values of holomorphic functions on the unit ball, continuous up to the boundary, and likewise for all the other spaces $\H_0^p,\,\P^p,\,\P_0^p,...$ etc. The Cauchy-Szego projection  sends $L^p$
into $\H^p$ boundedly.
\bigskip
\noin{\sl Sobolev spaces.}\smallskip
The Sobolev, or Folland-Stein, spaces on $\Hn$ and $S^{2n+1}$ can be defined in terms of the powers of the corresponding conformal sublaplacians.  The main references here are for example [ACDB], [ADB],  [Fo2]. We summarize the main properties below.

It is well known (see e.g. [St]) that for $Y_{jk}\in\H_{jk}$

$$\D Y_{jk}=\lambda_j\lambda_k Y_{jk},\qquad \lambda_j=j+{n\over2}\eqno(1.11)$$

For $F\in L^2(\Sn)$, we can write  $F=\sum_{j,k\ge0}\sum_{\ell=1}^{m_{jk}} c_{jk}^\ell(F) Y_{jk}^\ell$,  and $c_{jk}^\ell(F)=\ds\int F Y_{jk}^\ell$;  in particular, if  $F\in C^\infty(\Sn)$ then (1.11) implies that 
 $$\sum_{j,k\ge0}\sum_{\ell=1}^{m_{jk}}(\lambda_j\lambda_k)^d|c_{jk}^\ell(F)|^2<\infty.\eqno(1.12)$$
For $F\in C^\infty(\Sn)$ we then define for any $d\in\R$
$$\D^{d/2}F=\sum_{j,k\ge0}\sum_{\ell=1}^{m_{jk}}(\lambda_j\lambda_k)^{d/2} c_{jk}^\ell(F)Y_{jk}^\ell\eqno(1.13)$$
so that $\D^{d/2}$ extends naturally to the space of distributions on the sphere.
For $d>0,\,p\ge1$ we let 
$$W^{d,p}=\{F\in L^p: \,\D^{d/2}F\in L^p\}$$
endowed with norm 
$$\|F\|_{W^{d,p}}^{}=\|\D^{d/2}F\|_p;$$
the space $W^{d,p}$ is the completion of $C^\infty(\Sn)$ under such norm.

$W^{d,2}$ is the space of $F$ in $L^2$ so that (1.12) and (1.13) hold, and it is a Hilbert space 
with inner product and norm
$$(F,G)_{W^{d,2}}=\isn \D^{d/2}\!F\,\bar{\D^{d/2}G},\qquad \|F\|_{W^{d,2}}=(F,F)_{W^{d,2}}^{1/2}.$$
   Clearly $\|(I+\L)^{d/2}\|_2$ yields an equivalent norm on $W^{d,2}$.
Also, if $L_d^2$ denotes the classical Sobolev space on $\Sn$, defined as above but using the (positive) Laplace-Beltrami $\Delta$ rather than $\D$, and with norm $\|F\|_{L_d^2}=\|(I+\Delta)^{d/2}F\|_2$, then
$$L_d^2\hookrightarrow W^{d,2}\hookrightarrow L_{d/2}^2$$
in fact 
$$c_1\|F\|_{L_{d/2}^2}\le\|F\|_{W^{d,2}}\le c_2\|F\|_{L_d^2}$$
for some $c_1,c_2>0$, as one can easily see by comparing the eingenvalues
of $\D$ with those 
of $I+\Delta$ (i.e. $1+(j+k)(j+k+2n)$).

The dual of $W^{d,2}$ is  the space of distributions
$$(W^{d,2})'=\{\D^{d/2}F, F\in L^2\}$$
and it coincides with $W^{-d,2}$ defined as the space of distributions $T$ such that $\D^{-d/2}T\in L^2$.

The operators $\D^{d/2}$ and  $\L^{d/2}$ are positive and self-adjoint in their domain   $W^{d,2}$. The quadratic form $(\D^{d/4}F,\D^{d/4}G)$ allows us to further extend $\D^{d/2}$ and $\L^{d/2}$ to operators defined on $W^{d/2,2}$ (the form domain) valued in $W^{-d/2,2}$. In the sequel we will  denote such extensions by $\D^{d/2},\,\L^{d/2}$, with domain $W^{d/2,2}$.
\smallskip
On the Heisenberg group the Sobolev spaces are defined analogously as the completion of $C_c^\infty(\Hn)$ under the norm $\|(I+\L_0)^{d/2}\|_2$. The resulting space is  still denoted by $W^{d,2}$.
\bigskip

\noin{\sl Intertwining and Paneitz-type operators on the CR sphere}
\bigskip

The group $SU(n+1,1)$ acts  as a group of conformal transformations on $S^{2n+1}$,  and therefore on $\Hn$ by means of the Cayley projection (see [KR1-2]). Recall that a conformal (or contact)  transformation, is a diffeomorphism $h:\Hn\to\Hn$ that preserves the contact structure, i.e. if $\theta_0$ is a contact form, then $h^*\theta_0=|J_h|^{2/Q}\theta_0$, where $|J_h|$ is the Jacobian determinant of $h$. An analogue of the Euclidean Liouville's theorem holds: every $C^4$ conformal mapping on $\Hn$  comes from the action of an element of $SU(n+1,1)$, and it can be written as composition of 
$$\eqalign{{\hbox {left translations }} \;\;&(z,t)\to (z',t')(z,t)\cr
{\hbox {dilations }} \;\;&(z,t)\to (\delta z, \delta^2t),\,\;\delta>0\cr
{\hbox {rotations }} \;\;&(z,t)\to (Rz,t),\;\;R\in U(n)\cr
{\hbox {inversion }} \;\;&(z,t)\to \Big(-{z\over|z|^2+it},-{t\over |z|^4+t^2}\Big).\cr}$$
Let us denote  the spaces of conformal transformations (also called  CR automorphisms) of $\Hn$  by $\con(\Hn)$, and the space of conformal tranformations of $S^{2n+1}$  by  $\con(S^{2n+1}):=\{\tau: \tau=\Ca\circ h\circ \Ca^{-1}$ some $h\in \con(\Hn)\}.$
Note that the inversion on $\Hn$ corresponds to the antipodal map $\z\to-\z$ on $S^{2n+1}$.

The functions $|J_h|$ with $h\in\con(\Hn)$,  are obtained from $|J_\Ca|$ by left translations and dilations and can be written as (cf [JL2])
$$|J_h(u)|={C\over \big| \,|z|^2+it + 2 z\cdot w+\lambda\big|^Q},\qquad C>0,\, w\in \C^n,\, \lambda\in\C,\,\Re\lambda>|w|^2,\, u=(z,t)\in\Hn.$$
From this formula it follows  that the family of functions $|J_\tau|$ with $\tau\in \con(\Sn)$ can be parametrized as 
$$|J_\tau(\z)|={C\over |1-\omega\cdot\zeta|^Q},\quad C>0,\, \omega\in \C^{n+1},\, |\omega|<1,\, \z\in \Sn.\eqno(1.14)$$

The following formulas hold:
$$\eqalign{d(h(u),h(v))&=d(u,v)|J_h(u)|^{1\over2Q}|J_h(v)|^{1\over2Q},\qquad \forall h\in\con(\Hn)\cr 
d(\tau(\z),\tau(\eta))&=d(\z,\eta)|J_\tau(\z)|^{1\over2Q}|J_\tau(\eta)|^{1\over2Q},\qquad \forall \tau\in\con(\Sn)\cr}\eqno(1.15)$$

These formulas are trivially checked  on traslations, rotations, dilations of $\Hn$, and on the inversion of $\Sn$; using (1.2) one can  cover the remanining  cases. 

The operators $\L_0$ and $\D$ are intertwining in the sense that for each $f\in C_0^\infty(\Hn)$ and $F\in C^\infty(\Sn)$
$$|J_h|^{Q+2\over2Q} (\L_0 f)\circ h= \L_0\big(|J_h|^{Q-2\over2Q} (f\circ h)\big),\qquad \forall h\in\con(\Hn)$$
$$|J_\tau|^{Q+2\over2Q} (\D F)\circ \tau= \D\big(|J_\tau|^{Q-2\over2Q} (F\circ \tau)\big), \qquad\tau\in\con(\Sn).\eqno(1.16)$$

To check these formulas it is enough to rewrite them in terms of the inverse operators $\L_0^{-1},\,\D^{-1}$, and then use
the explicit formulas for their kernels and (1.15).

For $0<d<Q$ the general intertwining operator $\A_d$ of order $d$ is defined  by the following property:
$$|J_\tau|^{Q+d\over2Q} (\A_d F)\circ \tau= \A_d\big(|J_\tau|^{Q-d\over2Q} (F\circ \tau)\big),\;\;\forall\tau\in \con(S^{2n+1})\eqno(1.17)$$
for each $F\in C^\infty(\Sn)$. In other words, the  pullback of $\A_d$ by a conformal transformation $\tau$ satisfies
$$\tau^*\A_d (\tau^{-1})^*=|J_\tau|^{-{Q+d\over2Q}}\A_d |J_\tau|^{Q-d\over2Q}$$
where $\tau^*F=F\circ\tau$.

The concept of intertwining operator is more properly understood in the context of representation theory of semisimple Lie groups, in our case $SU(n+1,1)$, see e.g.  [Br], [BO\O], [C], [JW].  In particular,  for $d\in\C$ the map $u_d:\tau\to\big\{F\to |J_\tau|^{(Q+d)/(2Q)}(F\circ\tau)\big\}$ is a representation of the group $SU(n+1,1)$, modulo identification of the latter with $\con(\Sn)$; these $u_d$ are known as {\sl principal series representations} of $SU(n+1,1)$, and the ones corresponding to $d\in(-Q,Q)$ are called {\sl complementary series}. The relation (1.17) says that $\A_d$ intertwines the representations $u_{d}$ and $u_{-d}$. The present formulation is given in elementary differential-geometric terms, which for our purposes is more than enough (see however [Br], pp 18-19, for a digression on 
the $u_d$ in more Lie-theoretic language).  

It is known, from the above works (see also Appendix A, Prop. A.1), that an operator satisfying (1.17) is diagonal w.r. to the spherical harmonics, and its spectrum  is completely determined up to a multiplicative constant by the  functions
$$\lambda_j(d)={\Gamma\big({Q+d\over4}+j\big)\over\Gamma\big({Q-d\over4}+j\big)}\sim j^{d/2}\eqno(1.18)$$
in the sense that up to a constant the spectrum is precisely $\{\lambda_j(d)\lambda_k(d)\}$. From now on we will choose such constant to be 1, i.e. 
$\A_d$ will be  the operator on $W^{d,2}$ such that 
$$\A_d Y_{jk}=\lambda_j(d)\lambda_k(d)Y_{jk},\,\qquad Y_{jk}\in \H_{jk}\eqno(1.19)$$

The form  $(\A_d^{1/2}F,\A_d^{1/2}G)$ allows us to extend $\A_d$ to an operator with domain $W^{d/2,2}$ valued in $W^{-d/2,2}$, which we still denote by $\A_d$. The eigenvalues of such operators are still $\lambda_j(d)\lambda_k(d)$, i.e. (1.19) holds,   in the sense of forms. Since $\lambda_j(d)>0$ for all $j\ge0$ then $\ker \A_d=\{0\}$, and eigenvalue estimates show\ easily that   $\|\A_{d/2}F\|_2$ or $\|(\A_d)^{1/2}F\|_2$ are equivalent to $\|F\|_{W^{d,2}}$, for $0<d<Q$.
Observe that in the case $d=2$ we have $\lambda_j(2)=\lambda_j=j+{\ts{n\over2}}$, and we recover the conformal sublaplacian  i.e. $$\A_2=\D.$$

A fundamental solution of $\A_d$ is given by 
$$G_d(\z,\eta):=\A_d^{-1}(\z,\eta)=\sum_{j,k\ge0}
{\Phi_{jk}(\zbe)\over\lambda_j(d)\lambda_k(d)}=c_d \,d(\z,\eta)^{d-Q}\eqno(1.20)$$
with $$c_d={2^{n-{d\over2}}\,\Gamma\big({Q-d\over4}\big)^2\over \pi^{n+1} \Gamma\big({d\over2}\big)}\eqno(1.21)$$
and where the series converges unconditionally  in the sense of distributions, and also in $L^2$ if $Q/2<d<Q$.
The proof of (1.20) is somehow implicit in the work  of Johnson and Wallach [JW], and a similar formula (still quoted from [JW]) appears in  [ACDB] (formula (11)), but with different normalizations. The case $d$ an even integer was treated by Graham [Gr],  including the expression for the   fundamental solution. For the reader's sake  in Appendix A we offer a self-contained proof 
of the spectral characterization of intertwining operators, in the sense of (1.17), and of formula (1.20), using  only the explicit knowledge of the zonal harmonics and Schur's lemma. We note here (but see also Appendix A) that the intertwining property can be checked 
directly using (1.20) and formulas (1.15), after casting (1.17) in terms of the inverse $\A_d^{-1}$. 
\smallskip

\def\ss{\scriptstyle}
We shall be concerned with the intertwining, Paneitz-type  operators of order $Q$. Noticing that 
$$\lambda_0(d)={\Gamma\big({Q+d\over4}\big)\over\Gamma({Q-d\over4}\big)}\sim{Q-d\over4}\Gamma\big({Q\over2}\big),\quad d\to Q\eqno(1.22)$$
we easily obtain from (1.17) that the operator $\A_Q:W^{d,2}\to \P^\perp$ defined as
$$\A_QF:=\lim_{d\to Q}\A_dF\eqno(1.23)$$
the limit being in $L^2$, satisfies for $F\in W^{Q,2}$
$$|J_\tau| (\A_Q F)\circ \tau= \A_Q (F\circ \tau),\;\;\forall\tau\in \con(S^{2n+1})\eqno(1.24)$$
or
$$\tau^*\A_Q (\tau^{-1})^*=|J_\tau|^{-1}\A_Q. \eqno(1.25)$$

The operator $\A_Q$ can be extended via its quadratic form to an operator, still denoted by $\A_Q$, with domain $W^{Q/2,2}$, kernel $\ker \A_Q=\P$,   valued in $\big(W^{Q/2,2}\big)'=W^{-Q/2,2}$. The identity (1.24) is still valid for $F\in W^{Q/2,2}$ and 
$$\A_Q Y_{jk}=\lambda_j(Q)\lambda_k(Q)Y_{jk}=j(j+1)...(j+n)k(k+1)...(k+n)Y_{jk}.$$ 
Observe that $\|(I+\A_Q)^{1/2}F\|_2$ is equivalent to $\|F\|_{W^{Q/2,2}}$ on the space $W^{Q/2,2}\cap \P^\perp$.
\smallskip

In the case $d$ an even integer it is possible to write down a more explicit formula for $\A_d$ as a product
of Geller's type operators. In fact, we can recover the operators found by Graham in [Gr]:

\proclaim Proposition {1.1}. If $d\le Q$ is an even integer, then $\A_d$ is a differential operator and 
$$\A_d=\cases{\ds{\prod_{\ell=0}^{{\ts{d\over4}}-1}}\Big(\D-{(2\ell+1)^2\over4}+i(2\ell+1)\T\Big)
\Big(\D-{(2\ell+1)^2\over4}-i(2\ell+1)\,\T\Big) &  if $\,\;\ds{d\over4}\in\N$\cr
\D\,\ds{\prod_{\ell=0}^{\ts{d-2\over4}}}\Big(\D-\ell^2+2i\ell\T\Big)
\Big(\D-\ell^2-2i\ell\,\T\Big)&  if $\,\;\ds\;{d-2\over4}\in\N$.\cr}$$\par

\pf Proof. We have 
$$\lambda_j(d)=\prod_{\ell=0}^{{\ts{d\over2}}-1} \big(\lambda_j+\ell-\ts{d\over4}+{1\over2}\big)$$
from which we have  that (recall: $\lambda_j=j+\nhalf$)
$$\lambda_j(d)\lambda_k(d)=\cases{\ds{\prod_{\ell=0}^{{\ts{d\over4}}-1}} \big(\lambda_j^2-(\ell+\half)^2\big)\big(\lambda_k^2-(\ell+\half)^2\big) &  if $\,\;\ds{d\over4}\in\N$\cr
\lambda_j\lambda_k\,\ds{\prod_{\ell=0}^{\ts{d-2\over4}}} \big(\lambda_j^2-\ell^2\big)\big(\lambda_k^2-\ell^2\big) & if $\,\;\ds\;{d-2\over4}\in\N$.\cr}$$
The proof is completed noticing that 
$\T Y_{jk}={i\over2} \,(j-k)Y_{jk},$ for $Y_{jk}\in \H_{jk},$
and that 
$\big(\lambda_j^2-b^2\big)\big(\lambda_k^2-b^2\big)=\big(\lambda_j\lambda_k-b^2+b(j-k)\big)\big(\lambda_j\lambda_k-b^2-b(j-k)\big).$\endpf

Note in particular that when $d=4$

$$\A_4=\Big(\L+{n^2-1\over4}\Big)^2+\T^2.$$

Also, note that since $\T^2=-|\T|^2$ then one can isolate the highest order derivatives in the above expression, counting $T$ as an operator of order 2, and obtain 
$$\A_d=|2\T|^{d/2}\,{\Gamma\big(\L|2\T|^{-1}+{2+d\over4}\big)\over\Gamma\big(\L|2\T|^{-1}+{2-d\over4}\big)}+{\hbox{ lower order derivatives}}.\eqno(1.26)$$
The formula above needs of course to be suitably interpreted, as $\T$ is invertible only on the space $\bigoplus_{j\neq k}\H_{jk}$. For $d$ not an even integer, we speculate that there might still be a way to make sense out of (1.26), as the ``leading 
operator" appearing in that formula, has the same form as the intertwinor on the Heisenberg group (see (1.33)).\smallskip
\noin{\bf Remark.} It is possible to show that a fundamental solution for $\A_Q:\P^\perp\to\P^\perp$  is given by 
$$\A_Q^{-1}(\z,\eta)={2\over \omegan \Gamma\big({Q\over2}\big)^2}\,\log^2{d^2(\z,\eta)\over2}$$
(up to a CR-pluriharmonic function).
This calculation can be effected using the explicit formula for the fundamental solution of $\A_d$, and differentiating twice with respect to $d$ at $d=Q$ (note that the constant $c_d$ has a pole of order two ad $d=Q$).
\bigskip

\medskip
\noin {\sl Conditional intertwinors.}\smallskip
Of particular importance  for us is the existence of another intertwinor of order $Q$ defined on $\P$, which we call the {\sl conditional intertwinor}. This is defined  by its action on the spherical harmonics in the following way:
$$\A_Q'Y_{j0}=\lambda_j(Q)Y_{j0}=j(j+1)...(j+n)Y_{j0},\qquad \A'_Q Y_{0k}=\lambda_k(Q)Y_{0k}.\eqno(1.27)$$
Observe that $\|(I+\A_Q')^{1/2}F\|_2$ is equivalent to $\|F\|_{W^{Q/2,2}}$ on $W^{Q/2,2}\cap\P$, so that $\A_Q'$ can be extended in the usual way to $W^{Q/2,2}\cap \P$.\smallskip
We summarize the properties of $\A_Q'$ in the following proposition.

\proclaim Proposition {1.2}. The operator $\A_Q'$  defined as in (1.27) is   positive semidefinite,  self-adjoint  on $W^{Q/2,2}\cap\P$, and  $\ker \A_Q'=\H_{00}$. For each $F\in C^\infty(\Sn)\cap \P$ we have 
$$\A_Q'F=-{4\over \GammaQ}\,{\p\over\p d}\bigg|_{d=Q}(\A_dF)=\lim_{d\to  Q}{1\over\lambda_0(d)}\A_dF\eqno(1.28)$$
and for every $\tau\in\con(\Sn)$
$$|J_\tau|(\A_Q'F)\circ \tau=
\A_Q'(F\circ \tau)+{2\over Q\Gamma\big({Q\over2}\big)}\A_Q\big(\log|J_\tau|(F\circ \tau)\big).
\eqno(1.29)$$
Moreover, $\A_Q'$ is a differential operator with 
$$\A_Q'F=\prod_{\ell=0}^n\big(2|\T|+\ell\big)F=\prod_{\ell=0}^n\big({\ts{2\over n}}\L+\ell\big)F,\qquad \forall F\in C^\infty(\Sn)\cap\P\eqno(1.30)$$
and it is injective on  $\P_0$ with  fundamental solution
$$G_Q'(\z,\eta):=(\A_Q')^{-1}(\z,\eta)=-{2\over n!\,\omega_{2n+1}}\,\log{ d^2(\z,\eta)\over2}.\eqno(1.31)$$
\par

Note that (1.29) says that the intertwining property in the form (1.24) or (1.25) continues to hold for $\A_Q'$, but modulo distributions that annihilate $\P$ (or modulo functions in  $\P^\perp$, if $F\in W^{Q,2}$). Also, $\A_Q'$ is an intertwining operator if seen as an operator from $\P$ to $L^2/\P^\perp$. In particular, the representations intertwined by $\A_Q'$ are the standard shift $\tau\to\{F\to F\circ\tau\}$, on $\P$,  and $\tau\to\{\,[F]\to [\,(F\circ\tau) |J_\tau|\,]\,\}$ on $L^2/\P^\perp$.

\pf Proof. The eigenvalues  of $\A_Q'$ vanish when $j=0$ or $k=0$,  hence $\ker A_Q'=\H_{00}$ (the constants). The first identity follows easily from (1.22). To prove  (1.29), it is enough to take  the $d$-derivative at $Q$ of (1.17):
$$|J_h|(\A_Q'F)\circ \tau-{2\over Q\Gamma\big({Q\over2}\big)}|J_\tau|\log|J_\tau|(\A_QF)\circ \tau=
\A_Q'(F\circ \tau)+{2\over Q\Gamma\big({Q\over2}\big)}\A_Q\big(\log|J_\tau|(F\circ \tau)\big)
$$
for each $F\in C^\infty(\Sn)\cap \P$.  We can trivially check (1.30) when $F$ is a spherical harmonic. The last statement (1.31) follows from the formula
$$G_Q'(\z,\eta)=\sum_{j=1}^\infty{\Phi_{j0}(\zbe)+\Phi_{0j}(\zbe)\over\lambda_j(Q)}=2\Re\sum_{j=1}^\infty{\Phi_{j0}(\zbe)\over\lambda_j(Q)}={2\over\Gamma\big({Q\over2}\big)\omega_{2n+1}}\,\Re\sum_{j=1}^\infty{(\zbe)^j\over j}.$$
\endpf

\smallskip
\noin{\sl Intertwining operators on the Heisenberg group.}
\smallskip
For completeness we say a few words for the case of the intertwining operators on $\Hn$. \smallskip
We already know from (1.7) that there is a direct connection between $\A_2=\D$ and $\L_0$, via the Cayley transform. To find the analogue situation for $\A_d$ one basically has to find the operator on $\Hn$ with fundamental solution  $|u|^{d-Q}$, since this operator is easily checked to be intertwining. This has been done by Cowling [C] and the result can be formulated as follows. Consider the $U(n)-$spherical functions
$$\Phi_{\lambda,k}(z,t)=e^{i\lambda t-|\lambda|\,|z|^2} L_k^{n-1}(|\lambda|\,|z|^2),\qquad\lambda\neq0,\,k=0,1,2,...$$
where $L_k^{n-1}$ denote the classical Laguerre polynomial of degree $k$ and order $n-1$.
These are the eigenfunctions of the sublaplacian $\L_0$ and of $T=\partial_t$:
$$\L_0 \Phi_{\lambda,k}=|\lambda|(2k+n)\Phi_{\lambda,k},\qquad T\Phi_{\lambda,k} =i\lambda \Phi_{\lambda,k}.$$
On $\Hn$ there is a notion of ``group Fourier transform", which on radial functions (i.e. functions depending only on $|z|$ and $t$) takes the form
$$\widehat f(\lambda,k)=\int_{\sHn} \Phi_{\lambda,k}(z,t) f(z,t)\,du,\qquad f\in L^1(\Hn).$$
With this notation we have 
$$\widehat{\L_0 f}(\lambda,k)=|\lambda|(2k+n)\widehat f(\lambda,k),\qquad \widehat {T f}(\lambda,k)=-i\lambda \widehat f(\lambda,k).$$

In analogy with the sphere situation, one can show that up to a multiplicative constant  there is a unique operator $\L_d$ such that 
$$|J_h|^{Q+d\over2Q}(\L_d f)\circ h=\L_d\big(|J_h|^{Q-d)\over2Q}(f\circ h)\big),\qquad \forall h\in \con(\Hn)$$
for $f\in  C^\infty(\Hn)$, and such $\L_d$ is characterized by (under our choice of the constant) 
$$\widehat {\L_d f} (\lambda,k)=2^{d/2}|\lambda|^{d/2}\,{\Gamma\big(k+{Q+d\over4}\big)\over \Gamma\big(k+{Q-d\over4}\big)}\widehat f(\lambda,k)=2^{d/2}|\lambda|^{d/2}\lambda_k(d)\widehat f(\lambda,k),\eqno(1.32)$$
or, otherwise put,
$$\L_d=|2T|^{d/2}\,\,{\Gamma\big(\L_0|2T|^{-1}+{2+d\over4}\big)\over \Gamma\big(\L_0|2T|^{-1}+{2-d\over4}\big)}.\eqno(1.33)$$
With this particular choice of the multiplicative constant we have 

$$\L_2=\L_0,\qquad \L_4=\L_0^2+T^2=\L_0^2-|T|^2$$
$$\L_d\Big((2|J_\Ca|)^{Q-d\over2Q}(F\circ\Ca)\Big)=(2|J_\Ca|)^{Q+d\over2Q}(\A_d F)\circ\Ca$$
and a fundamental solution of $\L_d$ is
$$\L_d^{-1}(u,v)=C_d\,|v^{-1}u|^{d-Q},\,\qquad C_d={\ts {1\over2}}\,c_d={2^{n-{d\over2}-1}\,\Gamma\big({Q-d\over4}\big)^2\over \pi^{n+1} \Gamma\big({d\over2}\big)}.\eqno(1.34)$$
The proofs of these facts are more or less contained in [C, Thm. 8.1], which gives the  computation of  the group Fourier transform of $|u|^{d-Q}$. Note however, that our proof of the corresponding  facts on the sphere (Appendix A) can easily be adapted to this situation.

We remark here that in the case $d$ an even integer the operator $\L_d$ coincides with the operator found by Graham in [Gr].

\smallskip
The intertwinors at level $d=Q$ on $\Hn$ are obtained in the same manner as those for the sphere. There is the operator 
$$\L_Q=\lim_{d\to Q} \L_d$$ 
whose kernel is the  space of boundary  values of pluriharmonic functions on the Siegel domain  (modulo identification of its boundary with $\Hn$). In terms of $\A_Q$ we have 
$$\L_Q(F\circ\Ca)=2|J_\Ca|(\A_Q F)\circ\Ca.\eqno(1.35)$$
\smallskip
For the conditional intertwinor, we recall that $f$ is the boundary value of a holomorphic (resp. antiholomorphic) function on the Siegel domain if and only if $\widehat f(\lambda,k)=0$ if $k\neq0$ or $\lambda<0$ (resp.  $\lambda>0$).  So for $f$ 
a smooth CR-pluriharmonic function on $\Hn$ we can define, in analogy with $\A_Q'$ and via (1.32),
$$\L_Q'f={-{4\over \GammaQ}}{\partial\over\ \partial d}\bigg|_{d=Q} \L_d f=\lim_{d\to Q} {1\over\lambda_0(d) }\,\L_d f=|2T|^{Q/2}f.$$
With this definition we have for a smooth $F\in\P$
$$2|J_\Ca|(\A_Q'F)\circ\Ca=\L_Q'(F\circ\Ca)+{2\over Q\GammaQ}\L_Q\big(\log(2|J_\Ca|)(F\circ \Ca)\big)$$
which basically says that the conditional intertwinor on $\Sn$ is nothing but $|2T|^{Q/2}$ on the $\Hn$-pluriharmonic functions, ``lifted" from $\Hn$ to $\Sn$ via the Cayley map (note that the second term on the right is orthogonal the the pluriharmonics). Also, we have 
$$|J_h|(\L_Q'f)\circ h=\L_Q'(f\circ h)+{2\over Q\GammaQ}\L_Q\big(\log|J_h|(f\circ h)\big),\quad h\in\con(\Hn)$$
analogous to (1.29).
\bigskip
\noin{\sl Intertwining operators and change of metric.}
\bigskip

The sublaplacian and conformal sublaplacian can be defined intrinsecally
on any compact, strictly pseudoconvex CR manifold $M$, in terms of the contact form $\theta$; see e.g. [JL1], [St]. In particular, the conformal sublaplacian $\D_\theta$, corresponding to the contact form $\theta$, satisfies the simple transformation formula
$$\D_{W\theta}^{}=W^{-{Q+2\over4}}\D_\theta^{} W^{{Q-2\over4}},\eqno(1.36)$$
for any positive, smooth function $W$ on $M$, where $Q=2n+2$ and $2n+1$ is the dimension of the manifold. 

General intrinsic constructions of higher {\it integer} order  CR invariant operators have been established by works of
 Fefferman, Gover, Graham, Hirachi ([Hi], [FH], [GG]). A special but important case is the fourth order {\it CR Paneitz operator} $P$ in dimension $Q=4$, introduced in [Hi], which    satisfies 
$$P_{W\theta}=W^{-1}P_\theta.$$
The CR Paneitz operator was also recently studied in [CCY].
  
It is natural to speculate that a similar theory could be devised for the conditional intertwinors, acting on pluriharmonic functions, which we introduced here only in the standard structure of $S^{2n+1}$. Rather than attempting an intrinsic construction of such operators, we will present a natural extension of $\A_Q'$ from the standard contact form $\theta$ of $S^{2n+1}$ to a ``conformally changed" form $W\theta$, motivated by the intertwining property  (1.29). We will be interested in studying eigenvalues inequalities of such operators later on, as part of the proof of the Beckner-Onofri's inequality (0.8) (see Prop.~3.6).

In order to motivate our construction, which will be carried over the whole family of intertwinors $\A_d$, first observe that if $\theta$ is the standard form on $\Sn$ then (1.36) implies that 
$\D_{W\theta}$ is a positive and self-adjoint operator, densely defined on  $L^2(\Sn,W^{Q/2}d\zeta)$. By standard facts (which will be recalled below)
$\D_{W\theta}$ has eingevalues $0<\lambda_j(W)\uparrow\infty$, and by the intertwining property (1.16) (see proof of Prop. 1.3 below) such eigenvalues are invariant under the conformal action that preserves $L^{Q/2}$ norms: 
$$\lambda_j(W)=\lambda_j\big((W\circ\tau)|J_\tau|^{2/Q}\big).$$

We can now extend all this to the operators $\A_d$ and $\A_Q'$. For $0<W\in C^\infty(\Sn)$ and $0<d\le Q$, the $L^2$ Hermitian products 
$$(F,G)=\isn F\bar Gd\zeta ,\qquad (F,G)_W:=\isn F\bar G W^{Q/d}d\zeta $$ define equivalent norms on  $L^2$. It follows that
  $\P$  is a  closed subspace of $L^2$ under the product $(F,G)_W$, and there exists a corresponding orthogonal projection $\pi_W$:
 $$\pi_W:L^2\to\P$$

\smallskip
\proclaim Proposition {1.3}. Let $W\in C^\infty(\Sn)$, with $W>0$.  For $0<d\le Q$ the operator $$\A_d(W):=W^{-{Q+d\over2d}}\A_d W^{{Q-d\over2d}}$$ 
satisfies 
 $$\big(\A_d(W)F,G\big)_W=\big(\A_dF,G\big)\qquad F,G\in C^\infty (S^{2n+1})\eqno(1.37)$$
and it  can be extended  to a self-adjoint operator on  $W^{d/2,2}$, which is positive definite if $d<Q$, and positive semidefinite if $d=Q$, with $\ker \A_Q(W)=\P$. There is a sequence $\{\psi_j^W\}$ of real-valued  eigenfunctions of $\A_d(W)$ that form an orthonormal basis of $L^2$  with respect  to  $(F,G)_W$. 
\smallskip
 The operator $\A_d(W)$ and its eigenvalues $\big\{\lambda_j(W)\big\}_0^\infty$ are conformally invariant, in the sense that if $\tau\in \con(\Sn)$ and $W_\tau=(W\circ\tau)|J_\tau|^{d/Q}$ then
 $$\tau^*\A_d(W) (\tau^{-1})^*=\A_d(W_\tau) \eqno(1.38)$$and
$$\lambda_j(W)=\lambda_j(W_\tau),\qquad j\ge0.\eqno(1.39)$$
\smallskip
The operator
$$\A_Q'(W):=\pi_W^{}W^{-1}\A_Q'\eqno(1.40)$$
 satisfies    
$$\big(\A_Q'(W)F,G\big)_W=\big(\A_Q'F,G\big),\qquad F,G\in C^\infty (S^{2n+1})\cap\P\eqno(1.41)$$
and it  can be extended to a self-adjoint, positive semidefinite  operator on $W^{Q/2,2}\cap\P$, with $\ker \A_Q'(W)=\H_{00}$.
There is a sequence $\{\phi_j^W\}$ of real-valued eigenfunctions of $\A_Q'(W)$ that form an orthonormal basis of $\P$ with respect  to the product $(F,G)_W$. \smallskip
The operator  $\A_Q'(W)$ and its eigenvalues $\big\{\lambda_j'(W)\big\}_0^\infty$ are conformally  invariant, in the sense that   if  $\tau\in \con(\Sn)$ and $W_\tau=(W\circ\tau)|J_\tau|$ then
$$ \tau^*\A_Q'(W)(\tau^{-1})^*=\A_Q'(W_\tau).\eqno(1.42)$$and
$$\lambda_j'(W)=\lambda_j'(W_\tau),\quad j\ge0.\eqno(1.43)$$
\par
\pf Proof. This proposition follows in a more or less straightforward way from the standard spectral theory of forms and operators on Hilbert spaces (e.g. see [Sh, Thm. 7.7]). For $0<d<Q$ identity (1.37) is obvious, and  $\big(\A_d(W)^{1/2}F,\A_d(W)^{1/2}F\big)_W\ge c\|F\|_{W^{d/2,2}}$, some $c>0$, and  we can find  an o.n. basis of eigenfunctions of $\A_d(W)$ for $L^2$. Clearly, since $\A_Q$ is real, such eigenfunctions can be chosen to be real-valued. Identity (1.38) follows  from the intertwining property (1.17), and implies that   if $\lambda$ is an eigenvalue of
$\A_d(W_\tau)$ with eigenfunction $\psi$, then $\lambda$ is also an eigenvalue of $\A_d(W)$, with eigenfunction $\psi\circ\tau^{-1}$, which is (1.39). The proof for the case $d=Q$  is similar, by considering the positive operators $I+\A_Q(W)$ and $I+\A_Q'(W)$.  Identity (1.42) follows from
$$\big(\A_Q'(W)(G\circ \tau^{-1})\circ \tau,\phi\big)_{W_\tau}=\big(\A_Q' G,\phi\big)=\big(\A_Q'(W_\tau)G,\phi\big)_{W_\tau},\qquad G,\phi\in W^{Q/2,2}\cap\P$$
which in turn is a consequence of the intertwining property (1.29), i.e. 
$$|J_\tau|\big(\A_Q'(G\circ\tau^{-1})\big)\circ \tau=
\A_Q'G+H,\qquad G\in W^{Q/2,2}\cap\P$$
for some $H\in \P^\perp$, and the fact that $G\circ\tau^{-1}\in \P$,
 since the conformal transformations are restrictions of biholomorphic mappings on the unit ball.
\endpf

The naturality of the operators $\A_d(W)$ and $\A_Q'(W)$ is expressed by the intertwining relations in (1.38), (1.42). In the case $d$ integer the operators $\A_d(W)$ coincide with those obtained intrinsecally in [FH],  [GG], [Hi], within the class of contact forms $\{W\theta\}$.

\bigskip
\eject
\centerline{\bf{2.   Adams and Moser-Trudinger inequalities on the CR sphere}}
\bigskip
In this section we establish  new sharp Moser-Trudinger inequalities on  $\Sn$.  Two particular cases of such estimates will be needed in the next section, for the proof of Beckner-Onofri's inequality (see Propositions 3.3 and 3.4), but we believe that other cases are    of independent interest. 
The first special result we will need is a sharp inequality of type (0.9) for the operator $B_{Q/2}=(\A_Q')^{1/2}$:
$$\int_{\Sn}\exp\bigg[{\omegan(n+1)!\over2}\bigg({|F|\over\|(\A_Q')^{1/2}F\|_2}\bigg)^{2}\,\bigg]d\z\le C_0\eqno(2.1)$$
 for all $F\in W^{Q/2,2}\cap\R\P$, with zero mean; this is a key estimate in order to show that  the Beckner-Onofri functional (0.8) is bounded below. We will in fact establish a  version of (2.1) that  is valid for more general spectrally defined operators acting on pluriharmonic functions or on Hardy spaces, since its proof does not really require the specific structure of the operator $(\A_Q')^{1/2}$.

The second main result that we will need has to do with (0.9) for the operator $B_{Q/2}=\L^{Q/4}$. For technical reasons  we will in fact need to use the spectrally modified operator $L_\lambda^{Q/4}=({2\over n}\L)^{Q/4}\pi+\sqrt \lambda\,\L^{Q/4}\pi^\perp$ ($\lambda>0$), and the following estimate:
$$\int_{\Sn}\exp\bigg[{\omegan(n+1)!\over2\big(1+{k_n\over\lambda}\big)}\bigg({|F|\over\|L_\lambda^{Q/4}F\|_2}\bigg)^{2}\,\bigg]d\z\le C_0\eqno(2.2)$$
for all $F\in W^{Q/2,2}$ with zero mean, and some specific constant $k_n>0$ depending only on $n$. Such estimate will be needed to prove an Aubin's type inequality, for functions with vanishing center of  mass (Proposition 3.4). The above estimate (2.2) will be a special case of a more general sharp Moser-Trudinger inequality valid for arbitrary real powers less than $Q$ of the operator $a\L\pi+b\L\pi^\perp$ ($a,b>0$), which include the sublaplacian, in the same spirit as Adams' original results on $\Rn$ [Ad].

The main step in the proof of (2.1), (2.2), and their generalizations, is their equivalent formulation  in terms of suitable potentials, also known as  ``Adams' forms" of   Moser-Trudinger inequalities.

\bigskip

\noin{\sl Adams inequalities for convolution type operators on the CR sphere}\smallskip

Let us  introduce some  notation:
$$u=(z,t)\in \Hn,\quad\Sigma=\{u\in \Hn: |u|=1\} ,\quad u^*=(z^*,t^*)={u\over|u|}\in \Sigma$$
$$\z=\Ca(z,t)\in \Sn,\quad {1-\zn\over1+\zn}=|z|^2+it=|u|^2e^{i\theta}, \eqno(2.3)$$
$$\E= \Ca(\Sigma)=\big\{(\z_1,...,\z_{n+1})\in \Sn:\,\Re\z_{n+1}=0\big\}.$$
\smallskip
It is easy to see that a function $h(\z,\eta)$ is $U(n+1)-$invariant, i.e. $h(R\z,R\eta)=h(\z,\eta),\; \forall R\in U(n+1)$,
if and only if $h(\z,\eta)=g(\zbe)$ for some $g$ defined on the unit disk of $\C$. Furthermore, from (2.3) the function $g(\zbn)=g(\zn)$ is independent on $\Re \zn$, i.e. it is defined on $\E$, if and only if it is a function of the angle $\theta=\sin^{-1}t^*$.

A measurable function $\phi:\big[-{\pi\over2},{\pi\over2}\big]\to\R$ can be viewed as a function on $\Sigma$, via $\phi(\theta)=\phi(\sin^{-1}t^*)$, and we will use the notation
$$\int_\Sigma \phi\,du^*:=\int_\Sigma \phi(\sin^{-1}t^*)du^*=\omega_{2n-1}\int_{-\pi/2}^{\pi/2} \phi(\theta)(\cos\theta)^{n-1}d\theta\eqno(2.4)$$
whenever the integrals make sense. The formula on the right in (2.4) is easily checked via polar coordinates.
Finally, for $w\in\C,\, |w|<1$ we let  
$$\theta=\theta(w)=\arg{\,1-w\over1+w}\in\Big[-{\pi\over2},{\pi\over2}\Big]. $$

\medskip

\smallskip
\proclaim Theorem 2.1. Let $0<d<Q$ and $p=\ds{Q\over d}$.  Define
$$Tf(\z)=\isn G(\z,\eta)f(\eta)d\eta,\qquad f\in L^p(\Sn)$$
where 
$$\eqalign{G(\z,\eta)&=g\big(\theta(\zbe)\big)\,d(\z,\eta)^{d-Q}+O\big(d(\z,\eta)^{d-Q+\e}\big)=\cr&=
2^{{d-Q\over2}} g\big(\theta(\zbe)\big)\,|1-\zbe|^{d-Q\over2}+O\big(|1-\zbe|^{{d-Q+\e\over2}}\big),\quad\z\neq\eta\cr}\eqno(2.5)$$
for bounded and measurable $g:\big[\!-{\pi\over2},{\pi\over2}\big]\to\R$, with  
$\big|O\big(|1-\zbe|^{{d-Q+\e\over2}}\big)\big|\le C|1-\zbe|^{{d-Q+\e\over2}}$, some $\epsilon>0$, and with $C$ independent of $\z,\eta$.
\smallskip Then, there exists $C_0>0$ such that for all $f\in L^p(\Sn)$
$$\int_{\Sn}\exp\bigg[A_d\bigg({|Tf|\over\|f\|_p}\bigg)^{p'}\bigg]d\z\le C_0\eqno(2.6)$$
with
$$A_d={2Q\over\ds{\int_\Sigma}|g|^{p'} du^*}\eqno(2.7)$$
for every $f\in L^p(S^n)$, with ${1\over p}+{1\over p'}=1$. Moreover, if the function $g(\theta)$  is H\"olderian of order $\sigma\in(0,1]$ then the constant in (2.7) is sharp, in the sense that if it is replaced by a larger constant then there exists a sequence $f_m\in L^p(\Sn)$ such that the exponential integral in (2.6) diverges to $+\infty$ as $m\to\infty$. \par

 In [CoLu1] Cohn and Lu give a similar result in the context of the Heisenberg group, and for kernels of type $G(u)=g(u^*)|u|^{d-Q}$, i.e. without any perturbations. A version analogous to Theorem 2.1 can be stated and proved  also on  $\Hn$ (thus extending the result in [CoLu1]).

The point of Theorem 2.1. is that the expansion (2.5) is precisely that of the fundamental solutions of several (if not most) differential and pseudodifferental  operators of interest in CR geometry, including for example the sublaplacian and its powers. 
\smallskip

\pf Proof. 
The proof of this theorem is an application of  general results about Adams inequalities in measure-theoretic settings, recently obtained by Fontana and Morpurgo [FM]. In fact, (2.6) is an instant consequence  of Theorem 1 in [FM] and the following sharp asymptotic estimate on the distribution function of $G(\zeta,\eta)$:
 $$\big|\{\z:\,|G(\z,\eta)|>s\}\big| = s^{-{Q\over Q-d}}\,{1\over 2Q}\ints|g|^{Q\over Q-d}du^*+O\big(s^{-{Q\over Q-d}-\sigma}\big)\eqno(2.8)$$
 for a suitable $\sigma>0$, as $s\to+\infty.$ The proof of (2.8) is a ``routine" calculation based on the asymptotic expansion (2.5): first use the Cayley transform to reduce things to $\Hn$, then use polar coordinates to complete the job (see [BFM], Lemma 2.3 for details). 

The sharpness statement is proved in [FM], and follows  the same general philosophy
originally used by Adams and later by Fontana, Cohn-Lu and many others. In our case it is possible to check that the sequence $f_m$ in the statement of Theorem 2.1.  can be chosen as 
$$\!\!f_m(\eta)=\cases{ |G(\n,\eta)|^{d/(Q-d)}\,\sgn\big( G(\n,\eta)\big) & if $|G(\n,\eta)|\le m\;, \;d(\n,\eta)\ge 2 m^{-2/(Q-d)}$\cr \cr0 & otherwise.\cr}$$

\endpf

\medskip\smallskip
\noin {\sl Moser-Trudinger inequalities for operators of $d-$type on Hardy spaces}
\medskip
For a given  $d>0$, we say that a densely defined and self-adjoint operator $P_d$ on $\H$  is of   {\sl $d-$type} if
$$P_d Y_{j0}=\mu_{j0}Y_{j0},\qquad \forall Y_{j0}\in\H_{j0}\eqno(2.9)$$
for a given sequence $\{\mu_{j0}\}$ such that for $j\to\infty$
$$0\le\mu_{00}\le\mu_{10}\le\mu_{20}\le...\qquad\mu_{j0}=j^{d/2}+a_1 j^{d/2-\e_1}+...+a_{m} j^{d/2-\e_{m}}+O(j^{d/2-\e_{m+1}})\eqno(2.10)$$
for some  $0<\e_1<\e_2<...<{\e_{m+1}}$ with ${\ts{Q-d\over2}}<\e_{m+1}$. From this condition it follows that $\ker(P_d)$ is finite dimensional, and that $P_d$ is a continuous operator from $W^{d,2}\cap\H$ to $\H$. More generally, one defines operators of $d-$type on $\H^p$ as densely defined operators satisfying (2.9) and (2.10). Note that by (2.10) the operator $P_d$ can be written on $C^\infty\cap \H^p$  as a finite sum of powers of the sublaplacian, up to a smoothing operator. $P_d$ is a continuous operator from $W^{d,p}\cap\H^p$ to $\H^p$ and invertible if restricted to $\ker(P_d)^\perp$ with
$$\ker(P_d)^\perp:=\Big\{F\in \H^p:\; \isn F\phi_k=0,\, k=1,...m\Big\}$$
and where $\phi_1,...,\phi_m$ denote a basis of $\ker(P_d)$, the null space of $P_d$.
Operators of $d-$type on $\bar \H^p$ and $\P^p$ are defined similary, and the spectrum of such operators is denoted by $\{\mu_{0j}\}$ and $\{\mu_{j0},\,\mu_{0j}\}$ respectively, where the $\mu$'s satisfy a condition of type (2.10).

Clearly the operators $(\A_Q')^\alpha$ are of $\alpha Q-$type, for $\alpha>0$.

\smallskip

\proclaim Theorem {2.2}. If $P_d$ is an operator of $d$-type on  $\P^p$, with $0<d<Q$, then there is $C_0>0$ such that for any $F\in W^{d,p}\cap\P^p\cap \ker(P_d)^\perp$  and with $\ds{p={Q\over d},\;\;{1\over p}+{1\over p'}=1} $ we have 
$$\int_{\Sn}\exp\bigg[A_d\bigg({|F|\over\|P_dF\|_p}\bigg)^{p'}\bigg]d\z\le C_0\eqno(2.11)$$
with 
$$A_d={ 2Q\over \ds\int_\Sigma |g_d|^{p'} du^*}\eqno(2.12)$$ 
and $$g_d(\theta)={2^{{Q-d\over2}+1}\,\Gamma\big({Q-d\over2}\big)\over\omegan\,n!}\,\cos\big({\ts{Q-d\over2}}\,\theta\big).\eqno(2.13)$$ In the special case $d=Q/2$ (i.e. $p=p'=2$)
$$A_{Q/2}={ \omega_{2n+1}(n+1)!\over2}=(n+1)\pi^{n+1}\eqno(2.14)$$
and this constant is sharp, i.e. it cannot be replaced by a larger constant in (2.11).
\smallskip
If $P_d$ is of $d-$type on  $\H^p$, then for any $F\in W^{d,p}\cap\H^p\cap \ker(P_d)^\perp$ both (2.11) and (2.12) hold with $g_d=\ds{{2^{Q-d\over2}\Gamma\big({Q-d\over2}\big)\over n!\omegan}}$. In the special case $d=Q/2$  we have $A_{Q/2}=\omega_{2n+1}(n+1)!=2(n+1)\pi^{n+1}$ and this constant is sharp.

\par
\noin {\bf Remark.} Inequality (2.1) is a special case of (2.11).\smallskip 
\pf Proof. If  $P_d$ is of $d$-type on $\H^p$, then it is invertible on $F\in W^{d,p}\cap\H^p\cap \ker(P_d)^\perp$, and has fundamental solution defined by the formula
$$P_d^{-1}(\z,\eta):=\lim_{R\to 1^-} \sum_{\mu_{j0}\neq0}{\Phi_{j0}(\zbe)\over\mu_{j0}}\,R^j$$
in the sense of distributions and  pointwise for $\z\neq\eta$. 

Using that 
$$\sum_{\mu_{j0}\neq0}{\Phi_{j0}\over\mu_{j0}}\,R^j ={1\over n!\omegan}\,
\sum_{j\ge j_0}{\Gamma\big(j+{Q\over2}\big)\over\Gamma(j+1)}{(R\,\zbe)^j\over\mu_{j0}}
$$
and using the hypotesis on the $\mu_{j0}$ it is straightforward to check that 
$$P_d^{-1}(\zeta,\eta)={\Gamma\big({Q-d\over2}\big)\over\omegan\,n!}\,
(1-\zbe)^{d-Q\over2}+O\big(|1-\zbe|^{{d-Q\over2}+\e}\big)\eqno(2.15)$$
for a suitable $\e>0$.

Likewise,  if $P_d$ is of $d$-type on $\P^p$, then it is invertible on $F\in W^{d,p}\cap\P^p\cap \ker(P_d)^\perp$, and has fundamental solution defined by the formula
$$P_d^{-1}(\z,\eta):=\lim_{R\to 1^-} \bigg\{\sum_{\mu_{j0}\neq0}{\Phi_{j0}(\zbe)\over\mu_{j0}}\,R^j+
\sum_{\mu_{0j}\neq0}{\Phi_{0j}(\zbe)\over\mu_{0j}}\,R^j\bigg\}$$
in the sense of distributions and pointwise for $\z\neq \eta$, and the following expansion holds:
$$\eqalign{P_d^{-1}(\zeta,\eta)&={2\Gamma\big({Q-d\over2}\big)\over\omegan\,n!}\,\Re(1-\zbe)^{d-Q\over2}+O\big(|1-\zbe|^{{d-Q\over2}+\e}\big)\cr&=
2^{{d-Q\over2}}g_d(\theta)\,|1-\zbe|^{d-Q\over2}+O\big(|1-\zbe|^{{d-Q\over2}+\e}\big)\cr}\eqno(2.16)$$for a suitable $\e>0$.
Note that $(1-\zbe)=|1-\zbe|\,e^{i\theta}+O(|1-\zbe|^2).$

 The proof of (2.11) now follows from Theorem 2.1, taking $T$ to be the integral operator with kernel $G(\zeta,\eta)=P_d^{-1}(\zeta,\eta)$, as in (2.15) and (2.16). In the case $d=Q/2$ the computation of $A_{Q/2}$ is based on (2.4)
and the formula
$$\int_0^{\pi/2} \cos^2\big({\ts{n+1\over2}}\theta\big)(\cos\theta)^{n-1}d\theta={1\over2}\int_0^{\pi/2}(\cos\theta)^{n-1}d\theta={\sqrt\pi\,\Gamma\big({n\over2}\big)\over4\,\Gamma\big({n+1\over2}\big)}$$
together with the duplication formula for the gamma function.

For the proof of the sharpness statements see [BFM].
\endpf

\noin{\sl Moser-Trudinger inequalities for powers of sublaplacians}\bigskip
In this section we obtain sharp Moser-Trudinger inequalities  for  $\L^{d/2}$, and more generally  for powers of operators of type $L_{a,b}:=a\L\pi+b\L\pi^\perp$, where $\pi^\perp:=I-\pi$ on $L^p$. As in the proof of Theorem 2.2, the main step is to give precise asymptotic estimates for the fundamental solution of such operators. 

The starting point is an explicit formula for the fundamental solution of the powers of the $\Hn$ sublaplacian:
$$\L_0^{-d/2}(u,0)={\ts{1\over2}}\,G_d(\theta)|u|^{d-Q}$$
$$G_d(\theta)={2^{n+1}\Gamma\big({Q-d\over2}\big)\over \pi^{n+1}\Gamma\big({d\over2}\big)}\,\Re\bigg\{e^{i{Q-d\over2}\theta}\int_0^\infty\Big({s\over 1-e^{-2s}}\Big)^{{d\over2}-1} {e^{-ns}\over(e^{2i\theta}+e^{-2s})^{Q-d\over2}}\,ds\bigg\} \eqno(2.17)$$
which was derived first by [BDR] in case $d$ an even integer, and later by [CT] for any $d<Q$ using the heat kernel approach. 

\smallskip The following result  yields  more  information on the function $G_d(\theta)$, and it  will be useful in the explicit  computation of the sharp  constants for the case $p=2$.

\proclaim Proposition 2.3. $G_d(\theta)$ has the following trigonometric expansion
$$G_d(\theta)=\sum_{k=0}^\infty {g_{k,d}(\theta)\over \lambda_k^{d/2}}\eqno(2.18)$$
where
$$g_{k,d}(\theta)={2^{{Q-d\over2}+1}\over \omega_{2n+1}n!}\sum_{\ell=0}^k {(-1)^\ell\Gamma(k-\ell+d/2-1)\Gamma(\ell+n-d/2+1)\over\Gamma(d/2-1)\Gamma(k-\ell+1)\Gamma(\ell+1)}\;\cos\big[\big(2\ell+\ts{Q-d\over2}\big)\theta\big]$$
if $d\neq2$, with the series converging in the sense of distributions,  and
$$g_{k,2}(\theta)={(-1)^k2^{n+1}\over \omegan n!}\cdot{\Gamma(k+n)\over\Gamma(k+1)}.$$
 Moreover,
$$\int_{\Sigma} g_{k,d} g_{j,Q-d} du^*={4\,\Gamma(k+n)\over\pi^{n+1}\Gamma(n)\Gamma(k+1)}\,\delta_{j,k}^{}.
\eqno(2.19)$$ \par
  Formula (2.18) appeared in [BDR], for the case $d$ an even integer, and can be shown in a similar way but using formula (2.17), and writing  
$(1-e^{-2s})^{d/2-1}$ and $(e^{2i\theta}+e^{-2s})^{-(Q-d)/2}$ as binomial series. The orthogonality relation (2.19) seems  to be new, and its proof is a brute calculation involving classical terminating Saalsch\"utzian hypergeometric series  (see  [BFM] for more details).

\proclaim Proposition 2.4. The fundamental solution of $\L^{d/2}$ ($\,0<d<Q$) satisfies
$$\eqalign{\L^{-d/2}(\z,\eta)&=G_d(\theta)d(\z,\eta)^{d-Q}+O\big(d(\z,\eta)^{d-Q+\e}\big)\cr&=2^{{d-Q\over2}}G_d(\theta)|1-\zbe|^{d-Q\over2}+O\big(|1-\zbe|^{{d-Q+\e\over2}}\big)\cr}\eqno(2.20)$$
with $G_d(\theta)$ as in (2.17).  More generally, if 
$L_{a,b}:=a\L\pi+b\L\pi^\perp $ with $a,b>0,$
 then $L_{a,b}^{d/2}$ is continuous on $W^{d,p}$, invertible on the subspace of functions with  zero mean, and its fundamental solution  satisfies
$$L_{a,b}^{-d/2}(\z,\eta)=2^{{d-Q\over2}}\bigg[{g_d(\theta)\over(an/2)^{d/2}}+{g_d^\perp(\theta)\over b^{d/2}}\bigg]\,|1-\zbe|^{d-Q\over2}+O\big(|1-\zbe|^{{d-Q\over2}+\e}\big)\eqno(2.21)$$
$$g_d^\perp(\theta)=G_d(\theta)-{g_d(\theta)\over(n/2)^{d/2}}\eqno(2.22)$$
for a suitable $\e>0$, and with $g_d(\theta)$ as in (2.13), and $G_d(\theta)$ as in (2.17).
\par

Note that  $g_d(\theta)(n/2)^{-d/2}$ is the first term in the expansion (2.18), so that the notation $g_d^\perp$ in (2.22) is justified.

The proof of (2.20) is relatively straightforward in the case of integer powers, i.e. when $d$ is even. The idea is that first one should consider $\D^{-d/2}$, where $\D$ is the conformal sublaplacian with the explicit fundamental solution as in (1.6). The fundamental solution of $\D^{-d/2}$ is then a multiple integral on products of spheres, which can be related to the fundamental solution of $\L_0^{-d/2}$ on $\Hn$ via the Cayley transform. The case of $d$ not an even integer is more involved and the authors were able to handle it by using path integration. For details see [BFM], Proposition 2.6 and Corollary 2.7. The proof of (2.21) follows at once from (2.20) and the fact that the operator $\big({2\over n}\L\big)^{d/2}\pi$  is of $d-$type, so Proposition 2.3 applies.  \smallskip
The following is now  an immediate consequence of the above results combined with Theorem 2.1:

\medskip

\proclaim Theorem 2.5.  Let  $L_{a,b}=a\L\pi+b \L\pi^\perp$ ($\,a,b>0$). Then there is $C_0>0$ so that for any $F\in W^{d,p}$ with zero mean  and with $\ds{p={Q\over d},\;\;{1\over p}+{1\over p'}=1} $
$$\int_{\Sn}\exp\bigg[  A_d(a,b)\bigg({|F|\over\|L_{a,b}^{d/2}F\|_p}\bigg)^{p'}\bigg]d\z\le C_0\eqno(2.23)$$
with 
$$ A_d(a,b)={ 2Q\over \ds\int_\Sigma \bigg|{g_d(\theta)\over(an/2)^{d/2}}+{g_d^\perp(\theta)\over b^{d/2}}\bigg|^{p'}\,du^*}$$ 
and the constant $A_d(a,b)$ is sharp. If $d=\ds{Q\over2}$, or $p=p'=2$
$$A_{Q/2}(a,b)={\omegan (n+1)!\over2{\ds{\bigg[\Big({2\over an}\Big)^{n+1}+{1\over b^{n+1}}\sum_{k=1}^\infty
{{k+n-1}\choose{n-1}}\big(k+{\ts{n\over2}}\big)^{-n-1}\bigg]}}}.\eqno(2.24)$$
\par
\smallskip Setting $a=b=1$ in the above theorem gives the following sharp Moser-Trudinger inequality for the powers of the sublaplacian:\smallskip  
\proclaim Corollary 2.6. There is $C_0>0$ so that for any $F\in W^{d,p}$ with zero mean  and with $\ds{p={Q\over d},\;\;{1\over p}+{1\over p'}=1} $
$$\int_{\Sn}\exp\bigg[  A_d\bigg({|F|\over\|\L^{d/2}F\|_p}\bigg)^{p'}\bigg]d\z\le C_0$$
with 
$$ A_d={ 2Q\over \ds\int_\Sigma |G_d(\theta)|^{p'}\,du^*}\eqno(2.25)$$ 
and the constant $A_d$ is sharp. If $d=\ds{Q\over2}$, or $p=p'=2$
$$A_{Q/2}={(n+1)(n-1)!\,\pi^{n+1}\over{\ds{\sum_{k=0}^\infty
{(k+n-1)!\over k!\,\big(k+{\ts{n\over2}}\big)^{n+1}}}}}.\eqno(2.26)$$
In particular,  
$$A_{Q/2}=\cases{ 4 & if $n=1$\cr 18\pi & if $n=2$ \cr 192\,\ds{\pi^2\over 12-\pi^2} & if $n=3.$\cr}$$\par
\bigskip
\noin {\bf Remarks.} \smallskip\noin{{1.}} Inequality (2.2) is a special case of (2.23).\smallskip\noin{{2.}} The constant in (2.26) can be computed in principle  for any given $n$, by using partial fractions and the values of the 
Hurwitz zeta function $\sum_0^\infty (k+a)^{-s}$, when $a=n/2$ and $s$ is even. \smallskip \noin{{3.}}  
Corollary 2.6 above holds also for $\D^{d/2}$ with the same constant as in (2.25) (and for all functions in $W^{d,p}$). The reason for this is that the expansion (2.20) also holds for the kernel  of $\D^{-d/2}$ (see [BFM]. Prop. 2.6)
\bigskip\eject
\centerline{\bf 3. Beckner-Onofri's inequality}

\bigskip

The goal of this section is to establish the sharp Beckner-Onofri inequality for real CR-pluriharmonic functions on the sphere:

\proclaim Theorem 3.1. For any $F\in W^{Q/2,2}\cap\RP$ we have the inequality
$${1\over2(n+1)!}\avg F\A_Q'F\,d\z+\avg F\,d\z-\log \avg e^F\,d\z\ge0.\eqno(3.1)$$
The inequality is invariant under the conformal group of $\Sn$, in the sense that the functional on the left hand side is invariant under the action $F\to F^\tau=F\circ\tau+\log|J_\tau|$, for $\tau\in\con(\Sn)$.
Equality in (3.1) holds if and only if $F=\log|J_\tau|$, for some $\tau\in \con(S^{2n+1})$.
\par

There is a corresponding version of (3.1) for general complex-valued CR-pluriharmonic functions $F$:
$${1\over2(n+1)!}\avg \bar F\A_Q'F\,d\z+\avg \Re F\,d\z-\log \avg e^{{\sRe F}}\,d\z\ge0.$$
but it is a trivial consequence of the real-valued case.\smallskip

As we mentioned in the introduction, the proof of this theorem is based on the original compactness argument given by Onofri
 in dimension 2, and later perfected and extended to any dimensions by Chang-Yang, to provide an alternative proof of  Beckner's result. \smallskip

Define once and for all 
\def\J{{\cal J}}
$$\J[F]={1\over2(n+1)!}\avg F\A_Q'F\,d\z+\avg F\,d\z-\log \avg e^F\,d\z,$$
for any $F\in W^{Q/2,2}\cap\RP.$\smallskip
We divide the proof in three main steps:
\smallskip
\item{I.} Conformal invariance of $\J$\smallskip
\item{II.} Existence of a minimum for $\J$
\smallskip
\item{III.} Characterization of the minimum.

\bigskip
\noin{\sl Step I: Conformal invariance of $\J.$}
\bigskip
\proclaim Proposition 3.2. The conformal action $F\to F^\tau=F\circ\tau +\log|J_\tau|$ preserves $\RP$ and $W^{Q/2,2}\cap \RP$. Moreover, such spaces are the minimal  closed subspaces of  $L^2(\Sn),\, W^{Q/2,2}$ respectively, which are invariant under the conformal action. Finally, $\J[F^\tau]=\J[F]$, for all $F\in W^{Q/2,2}\cap\RP$.\par

\pf Proof. Clearly  $F\circ\tau\in\RP$
 if $F\in\RP$, and likewise for  $W^{Q/2,2}\cap\RP$.
For $\tau$ conformal, using (1.14) we see that  that $\log|J_\tau|\in\RP$. Any subspace $M$ of $L^2$
invariant under the action must contain the orbit of the function 0, i.e. all functions of type $\log|J_\tau|$; thus (still from (1.14)) every function of type $C-Q\,\Re\log (1-\z\cdot\omega)$  must be in $M$, for any given $\omega\in \C^{n+1},\,|\omega|<1$.
If $M$ is also closed, then it contains all   $\omega$-partial derivatives of such functions, evaluated at $\omega=0$, and therefore $M$ contains every real pluriharmonic polynomial and hence all of $\R\P$.\

Next consider the functional
$$\J_d[G]={1\over\lambda_0(d)^2}\avg G\A_d G\,d\z-\bigg(\avg|G|^{1/\theta}d\z\bigg)^{2\theta}$$
with $\theta={\ds {Q-d\over 2Q}}$. This functional is invariant under the action $G\to G_{\tau,\theta}=(G\circ\tau)|J_\tau|^{\theta}$; this follows from (1.17). One easily checks that as $\theta\to0$ (i.e. $d\to Q$)
$$\J_d[1+\theta F]={\theta^2\over\lambda_0(d)^2}\avg F\A_dF\,d\z+2\theta\avg F\,d\z-2\theta\log\avg e^F d\z+O(\theta^2)$$
so that if $F\in W^{Q/2,2}\cap \RP$, using (1.28) we obtain
$${d\over d\theta}\bigg|_{\theta=0}\J_d[1+\theta F]=2\J[F]$$
On the other hand, letting $G=1+\theta F$ we get  $G_{\tau,\theta}=(1+\theta F)_{\tau,\theta}=1+\theta F^\tau+O(\theta^2)$ so that 
$\J_d[(1+\theta F)_{\tau,\theta}]=\J_d[1+\theta F^\tau]+O(\theta^2)$
and by differentiation this implies $\J[F]=\J[F^\tau]$, if $F\in W^{Q/2,2}\cap \RP$.\endpf
\noin{\bf Note.} On the Euclidean $S^n$ the minimal subspace of $L^2$ that is invariant under the conformal action is the whole $L^2$. Indeed, in that case,  the $\log |J_\tau|$ are of type $C-n\log|1-\omega\cdot\zeta|$, with $\omega\in\R^{n+1},\,|\omega|<1$. An argument similar to the one used in the above proof shows that the orbit of the function $0$ is dense in $L^2$.

We remark that the proof above is an adaptation of Beckner's argument in [Bec]. Another possible proof of Prop. 3.2 can be given directly as in [CY], without appealing to the intertwining property of $\A_d$, but working directly with $\A_Q'$. We chose Beckner's
argument since it shows how the putative sharp, conformally invariant  Sobolev  inequality $\J_d[G]\ge0$ i.e.
$$\avg G\A_d G\,d\z\ge\bigg[{\Gamma\big({Q+d\over4}\big)\over\Gamma({Q-d\over4}\big)}\bigg]^2\bigg(\avg |G|^q\,d\z\bigg)^{2/q},\; \qquad q={2Q\over Q-d}\eqno(3.2)$$
would imply Beckner-Onofri's inequality (3.1), for pluriharmonics functions. Inequality (3.2), or its dual ``Hardy-Littlewood-Sobolev" form, is only known for $d=2$ [JL1,2]$^\dagger$\vfootnote\dag {See the ``Addendum" at the end of the Introduction.}
\bigskip

\noin{\sl Step II: Existence of a minimum for $\J$.}\bigskip
From now one we will denote the average of $F\in L^1(\Sn)$ by 
$$\wtilde F=\avg F={1\over\omegan} \isn F.$$\bigskip
\proclaim Proposition 3.3 (Provisional Beckner-Onofri's inequalities). There exists a constant $C$ such that for all $F\in  W^{Q/2,2}\cap\RP$ 
$${1\over 2(n+1)!}\avg F\A_Q'F\,d\z+\avg F-\log\avg e^Fd\z+C\ge0.\eqno(3.3)$$
If $\lambda>0$  then there exists  a constant $C_\lambda$ such that for all  $F\in W^{Q/2,2}$ and with $L_\lambda={2\over n}\L\pi+\lambda^{2/Q}\L\pi^\perp$
$$A_n(\lambda)\avg FL_\lambda^{Q/2} F\,d\z+\avg F-\log\avg e^Fd\z+C_\lambda\ge0\eqno(3.4)$$
 with
$$A_n(\lambda)={1\over 2(n+1)!}\,\bigg[1+{1\over\lambda}\,\sum_{k=1}^\infty{(k+n-1)!\over (n-1)!k!\big(k+{n\over2}\big)^{n+1}}\bigg].\eqno(3.5)$$\par

\medskip
\pf Proof. This is a standard argument based on the Adams inequalities (2.11) and (2.23) for the operators $(\A_Q')^{1/2}$ and $L_\lambda^{Q/4}=({2\over n}\L)^{Q/4}\pi+\sqrt \lambda\,\L^{Q/4}\pi^\perp.$ If an inequality of type
$$\int_{\Sn}\exp\bigg(B\,{|F-\wtilde F|^2\over\|PF\|_2^2}\bigg)d\z\le C_0$$
holds for one of the above  operators $P$ and  for either $W^{Q/2,2}\cap \RP$ or $F\in W^{Q/2,2}$ and with zero mean, then  letting $\mu=B^{1/2}(F-\wtilde F),\;\nu={1\over2} B^{-1/2}\|PF\|_2^2$ and expanding  $(\mu-\nu)^2\ge0$ we get 
$${1\over 4B}\,\|PF\|_2^2-\log\int e^{F-\wtilde F}d\z+\log C_0\ge0$$
which implies (3.3) and (3.4).
\endpf
\bigskip
\noin{\bf Remark.} We note that (3.3) is valid with $P_{Q/2}^2$ in place of $\A_Q'$, where $P_{Q/2}$ is any operator as in Prop. 2.3, with $d=Q/2$ and with kernel $\H_{00}$ (i.e. the constants).

\medskip
From (3.3) we now know that $\J$ is a functional that is  bounded below on $W^{Q/2,2}\cap \RP$. The goal now is to show that the minimizing sequence is actually bounded on such space. The first key step is the following Aubin's type inequality, used in the Euclidean setting first by Onofri and Aubin and then by Chang-Yang:
  
\proclaim Proposition 3.4 (Aubin's type inequality). For given $\sigma>{1\over2}$, there exist constants $C_1(\sigma),\,C_2(\sigma)$ such that for any $W^{Q/2,2}\cap \RP$ with $\isn \z_j \,e^Fd\z=0$ for $j=1,2...,n+1$, the following estimate holds 
$${\sigma\over2(n+1)!}\avg F \A_Q' F\,d\z+\avg F\,d\z-\log\avg e^Fd\z+C_1(\sigma)\|\L^{Q-1\over4}F\|_2^2+C_2(\sigma)\ge0\eqno(3.6)$$
\par

The proof below is an adaptation of the one in [CY, Lemma 4.6] (see also [Au2, Thm. 6]). We present it here because in our case there is an added 
difficulty, namely that the localization argument (multiplication by cutoff functions) inherent in the proof  does not 
preserve the class $\P$.

\bigskip
\pf Proof. Assume for the moment that $F\in W^{Q/2,2}$, and WLOG assume that $\int_{\Sn}e^F=\omega_{2n+1}$. Cover $\Sn$ with $2(2n+2)=2Q$ congruent spherical caps, by considering a cube inscribed inside the sphere, with side $L=2/\sqrt{2n+2}$. By rotation we can assume that if  
$$\Omega_{\delta_1}^1=\{x\in\Sn:\;\delta_1\le x_{2n+2}\le 1\},\qquad \delta_1<{1\over\sqrt{2n+2}}$$
then $$\int_{\Omega_{\delta_1}^1} e^F\ge {\omega_{2n+1}\over 2Q}\eqno(3.7)$$
It is not hard  to show that using the hypothesis $\int_{\Sn} x_{2n+2} e^F=0$, if
$$\Omega_{\delta_2}^2=\{x\in \Sn:\;-1\le x_{2n+2}\le -\delta_2\},\qquad \delta_2<{\delta_1\over 4Q}$$
then 
$$\int_{\Omega_{\delta_2}^2} e^F\ge\delta_2\omega_{2n+1}.\eqno(3.8)$$

Let $\phi_1,\phi_2$ be  cutoff functions such that $0\le\phi_j\le 1$ and
$$\phi_j=\cases{1 & on $\Omega_{\delta_j}^j$\cr\cr 0 & on $\Sn\setminus\Omega_{\delta_j/2}^j$\cr}$$

Consider  the  operator $L_\lambda={2\over n}\L\pi+\lambda^{2/Q}\L\pi^\perp$, 
so that from (3.4), (3.7) we obtain
$$\eqalign{{\omega_{2n+1}\over 2Q} &\le\int_{\Omega_{\delta_1}^1} e^F\le e^{\wtilde F}\int_{\Omega_{\delta_1}^1} e^{(F-\wtilde F)\phi_1}\le e^{\wtilde F}\omega_{2n+1}\avg e^{(F-\wtilde F)\phi_1}\cr&\le \omega_{2n+1}e^{\wtilde F} e^{C_\lambda}\exp\bigg[A_n(\lambda)\avg
(F-\wtilde F)\phi_1 L_\lambda^{Q/2}(F-\wtilde F)\phi_1+\avg(F-\wtilde F)\phi_1\bigg]\cr}\eqno(3.9)$$
with $A_n(\lambda)$ as in (3.5), and likewise, using (3.4) and (3.8) 
$$\delta_2\omega_{2n+1}\le\omega_{2n+1} e^{\wtilde F} e^{C_\lambda}\exp\bigg[A_n(\lambda)\avg
(F-\wtilde F)\phi_2 L_\lambda^{Q/2}(F-\wtilde F)\phi_2+\avg(F-\wtilde F)\phi_2\bigg].\eqno(3.10)$$
Now we claim that, for $k$ an even integer and $\e>0$
$$\eqalignno{\bigg|\isn(F-\wtilde F)\phi_j L_\lambda^k(F&-\wtilde F)\phi_j-{\big(\ts{2\over n}\big)^{k}} \isn \phi_j^2\big(\pi \L^{k/2}F \big)^2 -\lambda^{2k/Q}\isn \phi_j^2 \big(\pi^\perp \L^{k/2}F\big)^2\,\bigg|\cr&\le  \e\isn \big(L_\lambda^{k/2}F\big)^2+C(\lambda,\e)\isn F\L^{k-{1}} F,&(3.11)\cr}$$
whereas if $k$ is odd then 
$$\eqalignno{\bigg|&\isn(F-\wtilde F)\phi_j L_\lambda^k(F-\wtilde F)\phi_j-{\big(\ts{2\over n}\big)^{k}} \isn\phi_j^2 \big|\nabla_H\pi \L^{k-1\over2}F\big|^2-&(3.12)\cr&-\lambda^{2k/Q}\isn\phi_j^2 \big|\nabla_H\pi^\perp \L^{k-1\over2}F\big|^2 \bigg|\le  \e\isn \big(L_\lambda^{k/2}F\big)^2+C(\lambda,\e)\isn F\L^{k-{1}} F.\cr}$$
Here $\nabla_H$ denotes the so-called horizontal gradient defined on complex valued functions as 
$$\nabla_HF=\sum_{j=1}^{n+1}(\bar T_j\bar F\,T_j+T_j\bar F\,\bar T_j)$$
the $T_j$ being the generators  of $T_{1,0}(\Sn)$ defined in  (1.3). Such gradient satisfies the identities
$$\nabla_H G\cdot \bar{\nabla_H F}={1\over 2}\sum_{j=1}^{n+1}(\bar {T_jG}T_jF+T_j\bar G \,\bar T_j F)$$
$$\isn \bar G {\L F}=\isn \nabla_HG\cdot\bar {\nabla_H F}$$
Note that $\isn |\nabla_HL_\lambda^{k-1\over2} F|^2=\isn (L_\lambda^{k\over2} F)^2$.
The proof of these estimates  is given in the appendix, but the gist of it is that one can commute  $\phi_j$ with either the projection or $L_\lambda^k$, gaining one derivative of $F$. If $n$ is odd, using (3.11) (with $k=n+1$) 
we get for $j=0,1$
$$\eqalign{\isn(F-\wtilde F)\phi_j L_\lambda^{Q/2}(F-\wtilde F)\phi_j&\le \int_{\Omega_{\delta_j/2}^j}\bigg[{\big(\ts{2\over n}\big)^{k}} \big(\pi \L^{k/2}F \big)^2 +\lambda^{2k/Q} \big(\pi^\perp \L^{k/2}F\big)^2\bigg]\cr&+
\e\isn (L_\lambda^{Q\over4} F)^2+C(\lambda,\e)\|\L^{Q-1\over4}F\|_2^2.\cr}$$
Using these last inequalities in (3.9), (3.10) 
multipliying the resulting estimates out, and taking square roots we get
$$\sqrt{\delta_1\over 2Q}\le e^{\wtilde F}\exp \bigg[\Big(\half A_n(\lambda)+\e\Big)\avg FL_\lambda^{Q/2} F+C_1(\lambda,\e)\|\L^{Q-1\over4}F\|_2^2+C_2(\lambda)\bigg].$$
or
$$\Big(\half A_n(\lambda)+\e\Big)\avg FL_\lambda^{Q/2} F+\avg F+C_1(\lambda,\e)\|\L^{Q-1\over4}F\|_2^2+C_2(\lambda)\ge0$$
for some constants $C_1(\lambda,\e),C_2(\lambda)$. The case $n$ even is the same, just use (3.12) rather than (3.11).

Now, for  given $\sigma>{1\over2}$ we can certainly find $\lambda,\e$ so that $\ds{\half A_n(\lambda)+\e={\sigma\over2(n+1)!}}$, and  specializing  to $F\in W^{Q/2,2}\cap \RP$ we obtain 
$${\sigma\over2(n+1)!}\avg F \big({\ts{2\over n}}\L\big)^{Q/2} F\,d\z+\avg F\,d\z+C_1(\sigma)\|\L^{Q-1\over4}F\|_2^2+C_2(\sigma)\ge0.$$
Since on $\P$ we have $\big({\ts{2\over n}}\L\big)^{Q/2} \le \A_Q'$ we also obtain (3.6), under the condition $\ds{\avg e^F}=1$ (for the unconstrained case just replace $F$ in the above inequality by $F-\log{\ds\avg }e^F$). 

\endpf

We would like to make an important remark at this point. The very nature of the center of mass hypothesis in the above lemma makes it almost impossible to avoid the use of cutoff functions, in order to proceed with the localization argument; the authors were unable to conceive a different argument working exclusively inside the class $\P$.
This justifies our choice of the operator $L_\lambda$, which allows us to  temporarily exit the space $\P$.
Our choice is not the only one. For example, in the same spirit as in [CY] one could try to use the operator ${2\over n}\L$, i.e.   $L_\lambda$ with $\lambda^{4/Q}={2\over n}$. This operator satisfies  $\int F\big({2\over n}\L\big)^{Q/2} F\le \int F\A_Q'F$ for $F$ pluriharmonic, 
however to make the argument work the Adams constant $\wtilde A_{Q/2}$ corresponding to $\big({2\over n}\L\big)^{Q/2} $,  should satisfy $2\wtilde A_{Q/2}>A_{Q/2}$ with $A_{Q/2}$ as in (2.14). Using (2.24) we obtain 
 $${A_{Q/2}\over\wtilde A_{Q/2}}={\big({{n\over2}}\big)^{n+1}\over (n-1)!} \sum_{k=0}^\infty {(k+n-1)!\over k!\,\big(k+{n\over2}\big)^{n+1}}$$
which is less than 2 only for $n=1,2$ (in which cases one can indeed  use ${2\over n}\L$ to prove (3.6)), and seems to have exponential growth in $n$.

\medskip
The proof of the existence of the minimum for $\J$  can now proceed in more or less the same way as in [CY]. Let 
$$\S_0=\big\{F\in W^{Q/2,2}\cap \RP:\; \avg e^F\,d\z=1,\;\;\avg\z e^F\,d\z=0\big\}$$
and let us prove that a minimum of $\J$ exists  in $\S_0$.
First, we invoke the following  version of the ``center of mass theorem" for the CR sphere: if  $\ds\avg e^F=1$ then there exists  a conformal transformation $\tau$ such that  
$$\avg_{\Sn}\zeta \, e^{F^\tau(\zeta)}d\zeta =0,\qquad F^\tau=(F\circ\tau)+\log|J_\tau|.\eqno(3.13)$$
The proof of this fact is, by now,  a routine topological argument, modeled exactly after the proofs given in [CY1], [O], in the Euclidean case. The basic idea is that if the vector-valued integral (3.13) never vanishes as a function of $\tau$, then its unit normalization restricted to a suitable set of transformations can be seen as a   retraction from the closed unit ball of $\C^{n+1}$ to its boundary, which is not possible.

The  center of mass condition and the conformal invariance of $\J$  imply that  minimizing $\J$ over $W^{Q/2,2}\cap \RP$ is equivalent to minimizing $\J$ over $\S_0$.

Pick a minimizing sequence $F_k\in\S_0$, with $\J[F_k]\to\inf J$. Let's first prove that 
$$\avg F_k \A_Q' F_k\le C_2+C_1\|\L^{Q-1\over4}F_k\|_2^2.\eqno(3.14)$$
From (3.6), for a fixed ${1\over2}<\sigma<1$,
$$\J[F_k]+C_1(\sigma)\|\L^{Q-1\over4}F_k\|_2^2+C_2(\sigma)\ge{1-\sigma\over2(n+1)!}\avg F_k\A_Q'F_k$$
and since $F_k$ is minimizing we obtain (3.14). Now let's prove that $F_k$ can be chosen so that 
$$\|\L^{Q-1\over4}F_k\|_2\le C.\eqno(3.15)$$
For this we use the Ekeland principle (see e.g. [DeF], Thm 4.4.) to ensure that $\J'[F_k]\to0$ in $W^{-Q/2,2}\cap\RP$, where $\J'$ denotes the Gateaux  derivative of $\J$.  Thus,  $<\J'[F_k],\phi>=\int H_k\phi$ with 
$$H_k:=\A_Q'F_k-(n+1)!\,\pi(e^{F_k}-1)\to0 \;\;{\hbox{ in}}\;\; W^{-Q/2,2}\cap\RP$$ i.e.
$$F_k-\wtilde F_k=(\A_Q')^{-1}H_k+(n+1)! (\A_Q')^{-1}\pi(e^{F_k}-1)\eqno(3.16)$$

If $0<2\alpha<Q$, such as $\alpha={Q-1\over2}$, the operator $\A_Q'\L^{-\alpha/2}\pi$, with eigenvalues $\big({n\over 2}k\big)^{-\alpha/2}\lambda_k(Q)$, is of the type described by (2.9), (2.10), with $d=Q-\alpha$, hence  by Proposition 2.4 we have 
$$|\L^{\alpha/2}(\A_Q')^{-1}\pi (\z,\eta)|\le C|1-\zbe|^{-\alpha/2}.$$

So 
$$\eqalign{\isn &\big|\L^{\alpha/2}(\A_Q')^{-1}\,\pi (e^{F_k}-1)\big|^2d\z\le C\isn \bigg(\isn|1-\zbe|^{-\alpha/2}|e^{F_k(\eta)}-1|d\eta\bigg)^2d\z \cr&\le C\bigg(\isn\isn|e^{F_k(\eta)}-1|d\eta d\z\bigg)\isn\isn |1-\zbe|^{-\alpha}|e^{F_k(\eta)}-1|d\eta d\zeta \le C\cr} $$
(here we used that $\ds{\avg e^{F_k}=1}$ and that $\ds{\int|1-\zbe|^{-\alpha}}=C_\alpha$ for any $\eta\in \Sn$, since $2\alpha<Q$). On the other hand, looking at the eigenvalues of $\L^{\alpha/2}(\A_Q')^{-1}$
$$\isn \big|\L^{\alpha/2}(\A_Q')^{-1}H_k\big|^2d\z\le C\|H_k\|_{W^{\alpha-Q,2}}^2\le C\|H_k\|_{W^{-Q/2,2}}^2\le C$$
since $\|H_k\|_{W^{-Q/2,2}}\to 0$. All this with (3.16), $2\alpha=Q-1$,  and $\L^{\alpha/2}(F_k-\wtilde F_k)=\L^{\alpha/2}F_k$,  proves (3.15).

Finally, by Jensen's inequality $\wtilde F_k\le0$ and since $\J[F_k]\to\inf \J$ then
$$|\wtilde F_k|=-\avg F_k=-\J[F_k]+{1\over2(n+1)!} \avg F_k\A_Q' F_k \le C+ {1\over2(n+1)!} \avg F_k\A_Q' F_k\le C'$$
by (3.14) and (3.15). From this we deduce 
$$\avg |F_k|^2=\avg|F_k-\wtilde F_k|^2+|\wtilde F_k|^2\le C_1\|\L^{Q/4}F_k\|_2^2+C_2\le C$$
and therefore the minimizing sequence is bounded in $W^{Q/2,2}$. Now the standard argument goes like this: find a subsequence $F_{k_i}$ converging in $L^2$ and pointwise a.e. to an $F_0$, and weakly in $W^{Q/2,2}$. Clearly $F_0\in \RP$, and from the Adams inequality as $i\to\infty$, perhaps along another subsequence,
$$1=\avg e^{F_{k_i}}\to\avg e^{F_0}\qquad 0=\avg \z_j e^{F_{k_i}}\to \avg \z_j e^{F_{0}},\;\;j=1,2...,n+1$$
(this is because $e^{F_{k_i}}$ is  bounded in $L^2$, hence up to a subsequence it is weakly convergent, and its weak limit coincides with $e^{F_0}$ a.e.). Since $\ds\avg F_k\to\ds\avg F_0$ and $\J[F_k]$ converges, then also $\ds\avg F_k \A_Q' F_k$ converges, and  by standard results its limit is $\ge \ds  \avg F_0 \A_Q' F_0$, but it cannot be greater, since the sequence is minimizing for $\J$. This shows that $\J[F_k]\to\J[F_0]=\inf\J.$
\bigskip\eject
\noin{\sl Step III: Characterization of the minimum.}\bigskip

As in [CY] the problem of describing the minimum will be related to the first nonzero eigenvalue of a conformally invariant operator in the conformal class of the standard contact form, specifically the operator $\A_Q'(W)$ introduced in Prop. 1.3. 
According to Prop. 1.3, if $W\in C^\infty (\Sn)$ and $W>0$, then $\A_Q'(W)$ acting on  $W^{Q/2,2}\cap\P_0$, with inner product 
$(F,G)_W=\int FGW$, has positive eigenvalues $0<\lambda_1'(W)\le \lambda_2'(W)\le...$ (each counted with its multiplicity), and
$$\lambda_1'(W)=\min\bigg\{{(\phi,\A_Q'\phi)\over (\phi,\phi)_W^{}},\;\;\phi\in W^{Q/2,2}\cap\R\P_0,\, \isn \phi Wd\z=0\bigg\}\eqno(3.17)$$
Recall that $(\phi,\A_Q'\phi)=(\phi,\A_Q'(W)\phi)_W^{}$, and that the eigenfunctions of $\A_Q'(W)$ can be chosen real-valued.
\proclaim Proposition 3.5. Suppose that $F_0\in\S_0$ is a minimum for $\J$, 
then $F_0\in C^\infty(\Sn)$ and $\lambda_1'(e^{F_0})\ge (n+1)!$.\par

\pf Proof. The function $F_0$ must satisfy
$${d\over dt}\bigg|_{t=0}\J[F_0+t\phi]=\avg\phi\bigg({1\over2(n+1)!}\A_Q' F_0+ (e^{F_0}-1)\bigg)=0\qquad \forall\phi \in W^{Q/2,2}\cap\RP$$
i.e. ${1\over2(n+1)!}\A_Q' F_0+\pi (e^{F_0}-1)=0,$
which, together with (1.31),  implies that $F_0\in ~C^\infty(\Sn)$. On the other hand
$F_0$ must also satisfy
$${d^2\over dt^2}\bigg|_{t=0}\J[F_0+t\phi]={1\over(n+1)!}\avg \phi\A_Q'\phi+\bigg(\avg \phi e^{F_0}\bigg)^2-\avg \phi^2e^{F_0}\ge0$$
and from (3.17) we have $\lambda_1'(e^{F_0})\ge(n+1)!$.\endpf

The next result is a Hersch-type ``isoperimetric" inequality for the first $Q$ reciprocal eigenvalues. In the Euclidean case the inequality appeared first in [H] and it was later extended in [CY].
\smallskip
Notice that in our notation, when $W\equiv1$ on $\Sn$ we have 
$$\lambda_j'(1)=\lambda_1(Q)=(n+1)!,\;\;k=1,2,....,2n+2$$
since the bottom eigenvalue for $\A_Q'$ is $(n+1)!$ counted with multiplicity $m_{01}+m_{10}=2n+2$ (see (1.8)), its eigenspace being generated by the coordinate functions $\z_1,...,\zn$ and $\bar\z_1,...,\bar\zn.$

\bigskip\proclaim Proposition 3.6. For $W\in C^\infty(\Sn)$, $\,W>0$ and $\ds{\avg W=1}$  we have
$$\sum_{j=1}^{2n+2}{1\over\lambda_j'(W)}\ge \sum_{j=1}^{2n+2}{1\over\lambda_j'(1)}={2n+2\over\lambda_1(Q)}={2\over n!}\eqno(3.18)$$
In particular,
$$\lambda_1'(W)\le\lambda_1'(1)=(n+1)!\eqno(3.19)$$
and  equality holds in (3.18) or (3.19) if and only if $W=|J_\tau|$ for some $\tau\in \con(\Sn)$.\par
\bigskip\pf Proof. The proof of this uses the variational characterization of the sum of reciprocal eigenvalues (see [CY], or  [Ban, (3.7)])
$$\sum_{j=1}^{2n+2}{1\over\lambda_j'(W)}=\max\sum_{j=1}^{2n+2}{(\phi_j,\phi_j)_W^{}\over\big(\phi_j,\A_Q'(W)\phi_j\big)_W^{}}=\max\sum_{j=1}^{2n+2}{(\phi_j,\phi_j)_W^{}\over(\phi_j,\A_Q'\phi_j)}\eqno(3.20)$$
the maximum being over those $\phi_j\in W^{Q/2,2}\cap\P_0$ such that 
 $\;\ds{\avg \phi_jW=\avg\phi_j\bar{\A_Q'\phi_k}=0},$ for $j,k=1,...,2n+2,\;j\neq k.$ It is easy to see that the maximum is attained at $\phi_1,...,\phi_{2n+2}$ if and only if each $\phi_j$ is an eigenfunction of $\lambda_j'(W)$. 
By conformal invariance of the eigenvalues, i.e. $\lambda_j'(W)=\lambda_j'(W_\tau)$, where $W_\tau=(W\circ\tau)|J_\tau|$,  we can apply the center of mass theorem (3.13) with $W=e^F$, and  assume that $\ds{\avg \z_j W=0,}$ (and hence $\ds{\avg \bar\z_j W=0}$) for  $j=1,...,n+1$.

Therefore, we  can  choose $\z_j,\bar\z_j$ as
$\phi_j$ in (3.20), and since
$$(\z_j,\A_Q'\z_j)=\lambda_1(Q)\isn |\z_j|^2d\z={\omega_{2n+1}\over n+1}\,\lambda_1(Q)$$ we obtain
$$\sum_{j=1}^{2n+2}{1\over\lambda_j'(W)}\ge {n+1\over\lambda_1(Q)\omega_{2n+1}}\sum_{j=1}^{n+1}{\isn (|\z_j}|^2+|\bar\z_j|^2)W(\z)d\z={2(n+1)\over\lambda_1(Q)}$$
which is (3.18). Equality in (3.18) implies that each $\z_j,\;\bar\z_j$ is an eigenfunction of $\A_Q'(W)$ with eigenvalue $\lambda_1(Q)$, which implies  $\big(\phi,\A_Q'(W)\z_1\big)_W^{}=\lambda_1(Q)(\phi,\z_1)_W^{}$
for all $\phi\in C^\infty(\Sn)$, but this means $(\phi,\z_1)=(\phi,\z_1)_W^{}$ for all $\phi$ and this implies
$W\equiv1$ on $\Sn$. So if $W$ has vanishing center of mass then equality holds if and only if $W\equiv1$, so if we start from any $W$ by conformal invariance we have equality in (3.18) if and only if $W$ is in the conformal orbit of the constant function 1, i.e. $W=|J_\tau|$, some $\tau$.

Estimate (3.19) follows from the monotonicity of the eigenvalues, and equality in (3.19) implies equality in (3.18).
\endpf

To finish up the proof of Theorem 3.1, if $F_0\in \S_0$ is a minimum for $\J$ then $F_0\in C^\infty(\Sn)$
and by the previous propositions $\lambda_1'(e^{F_0})=\lambda_1(Q)=(n+1)!$, which is true if and only if 
$e^{F_0}=|J_\tau|$ for some $\tau\in \con(\Sn)$, and this concludes the proof.
\bigskip
\centerline{\bf 4. The logarithmic Hardy-Littlewood-Sobolev inequalities}
\medskip\medskip
In this final section we use the Beckner-Onofri inequality (3.1) to give a proof of (0.10), i.e. the CR version
of the inequality due to Carlen and Loss in the Euclidean setting [CL]. The procedure is at this point fairly standard, see for example [Bec] and [Ok]. The proof below is essentially the one in [Ok].

\bigskip


\proclaim Theorem {4.1} (Log HLS inequality). For any $G:\Sn\to \R$, with $G\ge0$, $G\in L\log L$, and ${\ds {\avg}} G=1$
we have 
$$(n+1)\avg\avg\log{1\over |1-\zbe|}\,G(\z) G(\eta)d\z d\eta\le \avg G\log G\,d\z\eqno(4.1)$$
with equality if and only if $G=|J_\tau|$,  for some $\tau\in\con (\Sn)$.\par

It is not hard to prove that for $G\in L\log L$ and $G\ge0$  the left-hand side is well-defined, finite and nonnegative. Also, 
in view of (1.31) when  $G\in L^2$   inequality (4.1) can be restated as
$${(n+1)!\over2}\avg (G-1) (\A_Q')^{-1}\pi (G-1)\le \avg G\log G.\eqno(4.2)$$


\medskip
Just like in  the Euclidean case it is possible to state an equivalent result on the Heisenberg group, via the Cayley transform:

\bigskip
\proclaim Corollary 4.2 (Log HLS inequality on $\Hn$). For any measurable $g:\Hn\to\R$ with $g\ge0$, $\ds{\int_{\sHn} g=\omega_{2n+1}}$ and $\ds{\int_{\sHn} g\log(1+|u|^2)<\infty}$ we have
$$(n+1)\avg_\sHn\avg_\sHn \log{2\over |v^{-1}u|^2} g(u)g(v)dudv\le \avg_\sHn g\log g+\log2\eqno(4.3)$$
where $\ds{\avg_\sHn={1\over\omega_{2n+1}}\int_{\sHn}}$. Equality in (4.3) occurs if and only if $g=(|J_\Ca|\circ h)|J_h|$  for some $h\in\con(\Hn)$.\par
\bigskip


\pf Proof of Theorem {4.1}. Let  $G\in L^2$, with $\, G\ge0$, $\,{\ds {\avg}} G=1$,
and let  $$F=(n+1)! \,(\A_Q')^{-1}\pi(G-1),$$ which is a real-valued pluriharmonic function with mean 0. Using Beckner-Onofri's inequality

$$\eqalign{{(n+1)!\over2}\avg (G-1)&(\A_Q')^{-1}\pi(G-1)={1\over2}\avg GF \cr& =\avg GF-{1\over2(n+1)!}\avg F\A_Q'F
\le\avg GF-\log\avg e^F.\cr}\eqno(4.4)$$
Now use Jensen's inequality to deduce
$$\log\avg e^F=\log \avg e^{F-\log G} G\ge\avg (F-\log G)G,\eqno(4.5)$$
which yields  (4.2)  for $G\in L^2$. Inequality (4.1) follows  for any  $G\in L\log L$ by an elementary truncation argument. From the Euler-Lagrange equation it is easy to see that any extremal of (4.1) must be in $C^\infty(\Sn)$.  Moreover, equality in (4.1) implies equality in (4.2), (4.4) and (4.5), i.e. (by Theorem 3.1) $F=\log|J_\tau|$ some $\tau\in \con(\Sn)$, and  $F-\log G=$constant, or $G=C|J_\tau|$; since $G$ has mean 1, then we finally have $G=|J_\tau|$ for some $\tau$. \endpf


\bigskip
\pf Proof of Corollary {4.2}. First observe that if $g:\Hn\to\R$ and $G:\Sn\to\R$ are related by $g=(G\circ\Ca)|J_\Ca|$ then 
 ${\ds\avg G=\avg_\sHn g=1}$ (with the above convention on the average on $\Hn$). Moreover, since $|1-\zbe|=2^{-{n\over n+1}}|J_\Ca(u)|^{1\over Q} |J_\Ca(v)|^{1\over Q}|v^{-1}u|^2$ (if $\Ca(u)=\z$ and  $\Ca(v)=\eta$) then
$$\eqalign{&(n+1)\avg\avg\log{1\over |1-\zbe|}\,G(\z) G(\eta)d\z d\eta- \avg G\log G\cr&=
(n+1)\avg_\sHn\avg_\sHn \log\Big(2^{{n\over n+1}}|v^{-1}u|^{-2}|J_\Ca(u)|^{-{1\over Q}}|J_\Ca(v)|^{-{1\over Q}}\Big)g(u)g(v)dudv-\cr&\hskip22em-\avg_\sHn g\log g+\avg_\sHn g\log|J_\Ca|\cr&=(n+1)\avg_\sHn\avg_\sHn \log{2\over |v^{-1}u|^2} g(u)g(v)dudv-\avg_\sHn g\log g-\log 2.
\cr}$$
This identity easily implies the statement.
The given integral condition on $g$ is to guarantee that $\int g\log g$ is finite if and only if $\int G\log G$ is finite, where $g$ and $G$ are related as above.\endpf

Note that with the same argument as in the proof of the Corollary above one can see that the log HLS functional (on $\Sn$ or $\Hn$)  is invariant under the conformal action.
\bigskip
\centerline{\bf 5. Appendix}
\bigskip\bigskip
\centerline{\bf A. Intertwining operators on $\Sn$}
\bigskip
In this appendix we give an explicit calculation of the spectrum of the intertwining operators $\A_d$, as defined by (1.17); a consequence of this calculation will be formula (1.20) up to a constant, and a further calculation will yield  the explicit constant given in  (1.21). The proof below is inspired by the method used by Johnson and Wallach [JW], 
but it is rather self-contained and uses a minimal apparatus from representation theory,  such as Schur's lemma and the knowledge of the zonal harmonics $\Phi_{jk}$. We believe that our calculation is actually sligthly simpler than that in [JW], at least in our context. In [Br] and [BO\O] there is another derivation of the spectrum of intertwining operators, again via the theory of spherical principal series representations of semisimple Lie groups ($SU(n+1,1)$ in our case), and the results there are quite general. 
\def\taul{{\tau_{\lambda}^{}}}
\proclaim Proposition A.1. Suppose that  an operator $\A_d$ $\,(0<d<Q)$ is intertwining, i.e. 
$$|J_\tau|^{Q+d\over2Q} (\A_d F)\circ \tau= \A_d\big(|J_\tau|^{Q-d\over2Q} (F\circ \tau)\big),\;\;\forall\tau\in \con(S^{2n+1})\eqno(5.1)$$
for each $F\in C^\infty(\Sn)$. 
Then $\A_d$ is diagonal with respect to the spherical harmonics, and  for every $Y_{jk}\in \H_{jk}$
$$\A_d Y_{jk}=c\lambda_j(d)\lambda_k(d)Y_{jk}$$
for some constant $c\in\R$, with $\lambda_j(d)$  as in (0.5). 
Viceversa, the operator $\A_d$ with eigenvalues $\lambda_j(d)\lambda_k(d)$ is intertwining, and has fundamental solution 
$$\A_d^{-1}(\z,\eta)=c_d d(\zeta,\eta)^{d-Q},\qquad c_d={2^{n-{d\over2}}\,\Gamma\big({Q-d\over4}\big)^2\over \pi^{n+1} \Gamma\big({d\over2}\big)}.$$\par

\pf Proof. For more clarity in the argument below we will use the notation $\H_{j,k},\,\Phi_{j,k}$ for $\H_{jk},\, \Phi_{jk},$ respectively. The fact that $\A_d$ is diagonal follows from Schur's lemma, and the irreducibility of the spaces $\H_{j,k}$.
 Suppose that $\A_d\Phi_{j,k}=\lambda_{j,k}\Phi_{j,k}$ with $\lambda_{j,k}=\lambda_{k,j}\in\R$
recall that $$\Phi_{j,k}(\z,\eta)=\Phi_{j,k}(\zbe):={(k+n-1)!(j+k+n)\over \omega_{2n+1}n!k!}(\bze)^{k-j}P_j^{(n-1,k-j)}(2|\zbe|^2-1)$$
if $j\le k$, and 
$\Phi_{j,k}(\z,\eta)=\Phi_{j,k}(\zbe):=\bar{\Phi_{k,j}(\zbe)}$, if $k\le j$. From now on we choose $\eta=\n$ and denote, for $j\le k$, 
$$\Psi_{j,k}(\zbn)=\Psi_{j,k}(z)=\bar z^{k-j}P_j^{(n-1,k-j)}(2|z|^2-1),\qquad z=\zbn=\zn$$
so that still $\A_d\Psi_{j,k}=\lambda_{j,k}\Psi_{j,k}$.

Consider the family of dilations of $\Hn$, which on the sphere take the form
$$\taul(\z)=\taul(\z',\zn)=\bigg({2\lambda\zeta'\over 1+\zeta_{n+1}+\lambda^2(1-\zeta_{n+1})},{ 1+\zeta_{n+1}-\lambda^2(1-\zeta_{n+1})\over  1+\zeta_{n+1}+\lambda^2(1-\zeta_{n+1})}\bigg)$$

The Jacobian of $\tau_\lambda$ is
$$|J_{\taul}|=\bigg|{2\lambda\over 1+z+\lambda^2(1-z)}\bigg|^Q$$
and 
$${d\over d\lambda}\bigg|_{\lambda=1}|J_{\taul}|^{a/Q}={a\over2}(z+\bar z).$$
Also, $\ds{{d\over d\lambda}\bigg|_{\lambda=1}(\taul \zbn)=z^2-1}$ so that  
$$\eqalign{{d\over d\lambda}\bigg|_{\lambda=1}|J_{\taul}|^{a/Q}&\big(\Psi_{j,k}\circ\taul\big)={a\over2}(z+\bar z)\bz^{k-j}P_j^{(n-1,k-j)}(2|z|^2-1)+\cr&+ (k-j)(-1+\bz^2)\bz^{k-j-1}P_j^{(n-1,k-j)}(2|z|^2-1)+\cr&+2(z+\bar z)(|z|^2-1)\bz^{k-j}{d\over dx} P_j^{(n-1,k-j)}(2|z|^2-1).\cr}\eqno(5.2)$$
The above quantity is a polynomial in $z,\bz$, with highest order monomials that are multiples of  $z^j\bz^{k+1}$ and $z^{j+1}\bz^{k}$. The projection of (5.2)
on $\H_{j+1,k}\bigoplus\H_{j,k+1}$ gives, for fixed $0\le j<k$
$$\eqalign{{d\over d\lambda}&\bigg|_{\lambda=1}|J_{\taul}|^{a/Q}\big(\Psi_{j,k}\circ\taul\big)\bigg|_{\H_{j+1,k}\oplus\H_{j,k+1}}=\cr&=A \bz^{k-j-1}P_{j+1}^{(n-1,k-j-1)}(2|z|^2-1)+B\bz^{k-j+1}P_j^{(n-1,k-j+1)}(2|z|^2-1)\cr}\eqno(5.3)$$
and for $j=k$
$$\eqalign{{d\over d\lambda}\bigg|_{\lambda=1}|J_{\taul}|^{a/Q}\big(\Psi_{j,k}\circ&\taul\big)\bigg|_{\H_{j+1,j}\oplus\H_{j,j+1}}=\cr&=A z P_{j}^{(n-1,1)}(2|z|^2-1)+B\bz P_j^{(n-1,1)}(2|z|^2-1)\cr}\eqno(5.4)$$
and the goal is to determine $A$ and $B$. In order to do this we consider the case $z$ real and $z$ imaginary, and compare the coefficients of the highest order powers in (5.2) and (5.3)-(5.4); what   we need here is that the coefficient of $x^j$ in a Jacobi polynomial of order $j$  is given by 
$${1\over j!}{d^j\over dx^j}P_j^{(\alpha,\beta)}(x)={1\over 2^j j!}{\Gamma(2j+\alpha+\beta+1)\over\Gamma(j+\alpha+\beta+1)}.$$
For $z$ real, a comparison of the coefficients of $z^{k+j+1}$ from $(5.2)$ and $(5.3)-(5.4)$ gives
$${\Gamma(k+j+n)\over j!\Gamma(k+n)}(a+k+j)=A\,{\Gamma(k+j+n+1)\over (j+1)!\Gamma(k+n)}+B\,{\Gamma(k+j+n+1)\over j!\Gamma(k+n+1)}$$
or
$$a+k+j=A\,{k+j+n\over j+1}+B\,{k+j+n\over k+n}.\eqno(5.5)$$
On the other hand, if $z$ is purely imaginary the same comparison yields
$$(-i)^{k-j+1}(k-j)\,{\Gamma(k+j+n)\over j!\Gamma(k+n)}=A(-i)^{k-j-1}{\Gamma(k+j+n+1)\over (j+1)!\Gamma(k+n)}+B(-i)^{k-j+1}\,{\Gamma(k+j+n+1)\over j!\Gamma(k+n+1)}$$
or
$$k-j=-A\,{k+j+n\over j+1}+B\,{k+j+n\over k+n}.\eqno(5.6)$$
Solving (5.5) and (5.6) in $A$ and $B$
$$A=\Big({a\over2}+j\Big){j+1\over k+j+n}\,,\qquad B=\Big({a\over2}+k\Big){k+n\over k+j+n}$$
which means, for $0\le j\le k$,
$$\eqalign{{d\over d\lambda}&\bigg|_{\lambda=1}|J_{\taul}|^{a/Q}\big(\Psi_{j,k}\circ\taul\big)\bigg|_{\H_{j+1,k}\oplus\H_{j,k+1}}=\cr&=\Big({a\over2}+j\Big){j+1\over k+j+n} \Psi_{j+1,k}+\Big({a\over2}+k\Big){k+n\over k+j+n}\Psi_{j,k+1}.\cr}\eqno(5.7)$$

Differentiating in $\lambda$ the intertwining relation (5.1) applied to $\Psi_{j,k}$ i.e. 
$$\lambda_{j,k}|J_\taul|^{Q+d\over2Q} (\Psi_{j,k}\circ \taul)= \A_d\big(|J_\taul|^{Q-d\over2Q} (\Psi_{j,k}\circ \taul)\big)$$
(it is easy to see that differentiation in $\lambda$ commutes with $\A_d$) and using (5.7) 
$$\eqalign{&\lambda_{j,k}\Big({Q+d\over4}+j\Big){j+1\over k+j+n} \Psi_{j+1,k}+\lambda_{j,k}\Big({Q+d\over4}+k\Big){k+n\over k+j+n}\Psi_{j,k+1}=\cr&=\lambda_{j+1,k}\Big({Q-d\over4}+j\Big){j+1\over k+j+n} \Psi_{j+1,k}+\lambda_{j,k+1}\Big({Q-d\over4}+k\Big){k+n\over k+j+n}\Psi_{j,k+1}\cr}$$
which implies
$$\lambda_{j+1,k}=\lambda_{j,k}\,{{Q+d\over4}+j\over{Q-d\over4}+j},\qquad \lambda_{j,k+1}=\lambda_{j,k}\,{{Q+d\over4}+k\over{Q-d\over4}+k}\;\quad k\ge j\ge0$$ 
and therefore
$$\lambda_{j,k}=\lambda_{0,k}\,{\Gamma\big({Q+d\over4}+j\big)\over\Gamma\big({Q-d\over4}+j\big)}=\lambda_{0,0}{\Gamma\big({Q+d\over4}+j\big)\over\Gamma\big({Q-d\over4}+j\big)}{\Gamma\big({Q+d\over4}+k\big)\over\Gamma\big({Q-d\over4}+k\big)}.$$

The proof of the last statement follows from the fact that the convolution operator $\B_d$ with kernel $d(\z,\eta)^{d-Q}$ is intertwining, but with $d$ replaced by $-d$: 
$$\B_d\big(|J_\tau|^{Q+d\over2Q}(G\circ\tau)\big)=|J_\tau|^{Q-d\over2Q}(\B_dG)\circ\tau$$
which can be checked directly on the dilations, translations (and trivially rotations and the inversion), using formulas (1.15).

From this and the previous calculations (which are valid also for  $-Q<d<0$) we deduce (note: $\lambda_j(-d)=\lambda_j(d)^{-1}$)
$$\isn d(\zeta,\eta)^{d-Q}Y_{jk}d\eta={c\over \lambda_j(d)\lambda_k(d)}Y_{jk}.$$
Now   set $j=k=0$, and by an  elementary computation 
$$\isn d(\z,\eta)^{d-Q}d\eta=2^{d-Q\over2}\isn|1-\zbe|^{d-Q\over2}d\eta=2^{{d-Q\over2}}\omega_{2n+1}\,{\Gamma\big({Q\over2}\big)\Gamma\big({d\over2}\big)\over\Gamma\big({Q+d\over4}\big)^2}$$
so that 
$$c=\lambda_0(d)^2\omega_{2n+1}\,{\Gamma\big({Q\over2}\big)\Gamma\big({d\over2}\big)\over2^{Q-d\over2}\Gamma\big({Q+d\over4}\big)^2}=\,\omega_{2n+1}\,{\Gamma\big({Q\over2}\big)\Gamma\big({d\over2}\big)\over2^{Q-d\over2}\Gamma\big({Q-d\over4}\big)^2}={1\over c_d}.$$
The operator $\A_d$ with eigenvalues $\lambda_j(d)\lambda_k(d)$ is the inverse of $c_d \B_d$ and so  it is also intertwining and has the requested   fundamental solution.
\endpf
\bigskip

\centerline{\bf B. Proofs of {(3.11)} and (3.12)}
\bigskip In the proofs of (3.11) and (3.12) we will assume WLOG that $F$ has zero mean, since the operators on the right-hand sides of such inequalities both annihilate the constants. To start with (3.11), we assume  $k$ even. We have $L_\lambda^k =\big({2\over n}\big)^k \pi\L^k +\lambda^{2k/Q}\pi^\perp \L^k$, and (for $\phi \in C^\infty$) 
$$\isn \phi F L_\lambda^k \phi F={\big(\ts{2\over n}\big)^{k}} \isn \big[\pi \L^{k/2} (\phi F)\big]^2+\lambda^{2k/Q}\isn \big[\pi^\perp \L^{k/2}(\phi F)\big]^2\eqno(5.8)$$
so let us first consider the first term.  Using  the definition of $\L$ we can write 
$$\L^{k/2}(\phi F)=\phi\L^{k/2} F+\sum_I \phi_I T_I F$$
where the sum is finite, over a suitable set of multiindeces $I=\{i_1,...,i_\ell\}$, $\ell<{k}$, and where $T_I=T'_{i_1}...T'_{i_\ell}$, the  $T'_j$ being either $T_j$ or $\bar T_j$, and $\phi_I$ a smooth function.  Apply $\pi$ to this formula and square it; the leading term is $(\pi\phi\L^{k/2}F)^2$, and the remainder terms 
are estimated using the following inequalities:

\item {i)} $\;\|\pi G\|_2\le \|G\|_2$\smallskip
\item{ii)} $\;\|T_I F\|_2\le C\|\L^{k-1\over2} F\|_2$, if $I$ has length $<k$
\smallskip
\item{iii)} $\;\|\pi\L^{k/2}F T_I F\|_1\le \e\|\pi\L^{k/2}F\|_2^2+ C(\e)\|\L^{k-1\over2} F\|_2^2.$
\smallskip
For ii) see for example [ADB], for an o.n.  base of $T_{1,0}$ rather than the $T_j$. Observe that ii) is also 
valid for $I$ empty, i.e. for $\|F\|_2$, since $F$ has zero mean.

Thus we are reduced to estimate the last two terms of the identity
$$\int \big[\pi(\phi \L^{k/2}F)\big]^2=\int \phi^2 (\pi\L^{k/2}F)^2+\int \Big([\pi,\phi]\L^{k/2}F\Big)^2 +2\int\Big([\pi,\phi]\L^{k/2}F\Big)\phi\pi\L^{k/2}F, $$
where $[\pi,\phi]=\pi\phi-\phi\pi$. In order to do this we just have to justify that if $T_j$ is as in (1.3)
then the operator $T=T_j[\pi,\phi]$ (and hence $[\pi,\phi]T_j$) is bounded on $L^2$. This is a  consequence of the famous $T1-$theorem
by David-Journe, in the context of spaces of homogeneous type (such as the CR sphere); see for example [DJS]. Indeed one can write down  explicitly the kernel of such operator, using the Cauchy-Szego kernel, and check that it is a Calderon-Zygmund kernel, with $T1=T^*1=0$. 

 This given, we can easily estimate the second and third term with  $\,\e\|\pi\L^{k/2}F\|_2^2+C(\e)\|\L^{k-1\over2}F\|_2^2$.
This takes care of the first term on the right-hand side of (5.8); to deal with the second term in (5.8), we argue exactly in the same manner. 
This shows (3.11) in case $k$ even.

 For $k$ odd, the proof of (3.12) is completely similar, except one has to start from $\int \pi\L^{k-1\over2}(F\phi)\pi\L^{k+1\over2}(F\phi)$. Using the same product rule as above and the commutator estimate, the leading term is given by 
$$\int \phi^2\pi\L^{k-1\over2}F\,\pi\L^{k+1\over2}F=\int \phi^2\big|\nabla_H\pi\L^{k-1\over2}F\big|^2+
\int \pi\L^{k-1\over2}(F\phi)\nabla_H\phi^2\nabla_H\pi\L^{k-1\over2}F$$
and it is easy to see that the second term is bounded above by 
$$\e\int\big|\nabla_H\pi\L^{k-1\over2}F\big|^2+C(\e)\|\L^{k-1\over2}F\|_2^2=\e\|\L^{k/2}\pi F\|_2^2+C(\e)\|\L^{k-1\over2}F\|_2^2.$$ 
The remainder terms are estimated similarly. \endpf

\vskip1.5em\eject \centerline{\bf References}\bigskip
\item{[Ad]} Adams, David R. {\sl
A sharp inequality of J. Moser for higher order derivatives},
Ann. of Math. {\bf128} (1988), no. 2, 385--398. 
\smallskip
\item{[ACDB]} Astengo F., Cowling M., Di Blasio B., {\sl The Cayley transform and uniformly bounded representations},
 J. Funct. Anal. {\bf213} (2004),  241-269.\smallskip
\item{[ADB]}  Astengo F.,  Di Blasio B., {\sl Sobolev spaces and the Cayley transform}, Proc. Amer. Math. Soc. {\bf134} (2006), 1319-1329.\smallskip 
\item{[Au1]} Aubin T., {\sl Probl\`emes isop\'erim\'etriques at espaces de Sobolev}, J. Differential Geometry {\bf11}
 (1976), 573-598.
\item{[Au2]}Aubin T.,
{\sl Meilleures constantes dans le th\'eor\`eme d'inclusion de Sobolev et un th\'eor\`eme de Fredholm non lin\'eaire pour la transformation conforme de la courbure scalaire}, J. Funct. Anal. {\bf32} (1979), 148-174. 
\smallskip\item{[BMT]} Balogh Z.M., Manfredi J.J., Tyson J.T., {\sl  Fundamental solution for the $Q$-Laplacian and sharp Moser-Trudinger inequality in Carnot groups}, J. Funct. Anal. {\bf204} (2003), 35-49.\smallskip
\item{[Ban]} Bandle C.,
{\sl Isoperimetric inequalities and applications}, 
Monographs and Studies in Mathematics, 7. Pitman (Advanced Publishing Program), Boston, Mass.-London, 1980.\smallskip
\item{[Bec]} Beckner W., {\sl Sharp Sobolev inequalities
on the sphere and the Moser-Trudin-\break ger inequality}, Ann. of
 Math. {\bf138} (1993), 213-242.\smallskip
\item{[BDR]} Benson C., Dooley A.H., Ratcliff G., {\sl Fundamental solutions for powers of the Heisenberg sub-Laplacian}, Illinois J. Math. {\bf 37} (1993), 455-476.\smallskip
\item{[Br]} Branson T.P., {\sl Sharp inequalities, the functional determinant 
and the complementary series},  Trans. Amer. Math. Soc. 
{\bf347} (1995), 3671-3742. 
\smallskip
\item{[Br1]} Branson T.P., Memo to No\"el Lohou\'e, 1999.\smallskip

\item{[BCY]} Branson T.P., Chang S-Y.A., Yang P., {\sl
 Estimates and extremals for zeta function determinants on
four-manifolds}, Commun. Math.  Phys. {\bf149} (1992),
241-262.\smallskip
\item{[BFM]} Branson T.P., Fontana L., Morpurgo C., {\sl Moser-Trudinger and Beckner-Onofri's inequalities on the CR sphere} (2007), arXiv:0712.3905v3.\smallskip
\item{[BO\O]} Branson, T.P., \'Olafsson G., \O rsted B., {\sl
Spectrum generating operators and intertwining operators for representations induced from a maximal parabolic subgroup},
J. Funct. Anal. {\bf135} (1996),  163-205. \smallskip
\item{[CL]} Carlen E., Loss M., {\sl Competing
symmetries, the logarithmic HLS inequality and Onofri's inequality on
$S^n$}, Geom. and Funct. Anal. {\bf2} (1992), 90-104.\smallskip
\item{[CT]} Chang D.-C., Tie J., {\sl Estimates for powers of sub-Laplacian on the non-isotropic Heisenberg group}, J. Geom. Anal. {\bf10} (2000),  653-678. \smallskip
\item{[CQ]} Chang S.-Y.A., Qing J., {\sl  The zeta                
   functional determinants on manifolds with boundary. II. Extremal
   metrics and compactness of isospectral set}, J. Funct. Anal. {\bf147}
   (1997), 363-399.\smallskip
\item{[CY]} Chang S.-Y.A., Yang P., {\sl Extremal             
   metrics of zeta function determinants on $4$-manifolds}, Ann. of Math.
   {\bf142} (1995),  171-212. \smallskip
\item{[CY1]} Chang S.-Y.A.-Yang  P.C., {\sl  Prescribing gaussian curvature
 in $S^2$}, Acta Math. {\bf159} (1987), 215-259.\smallskip
\item{[CCY]} Chanillo S., Chiu H.-L., Yang P.,  {\sl Embeddability for Three-Dimensional Cauchy-Riemann Manifolds and CR Yamabe Invariants}, to appear in Duke Math. J.\smallskip
\item{[CoLu1]} Cohn W.S., Lu G., {\sl Best constants for Moser-Trudinger inequalities on the Heisenberg group}, Indiana Univ. Math. J. {\bf50} (2001), 1567-1591. \smallskip
\item{[CoLu2]}  Cohn W.S., Lu G., {\sl Sharp constants for Moser-Trudinger inequalities on spheres in complex space $\C^n$}, Comm. Pure Appl. Math. {\bf57} (2004), 1458-1493. \smallskip
\item{[C]} Cowling M.,  {\sl Unitary and uniformly bounded representations of some simple Lie
groups},
   Harmonic Analysis and Group Representations, C.I.M.E.,
      Napoli: Liguori, (1982), 49-128.
\smallskip
\item{[DeF]} de Figueiredo D.G., {\sl
Lectures on the Ekeland variational principle with applications and detours},
Tata Institute of Fundamental Research Lectures on Mathematics and Physics, 81. Published for the Tata Institute of Fundamental Research, Bombay; by Springer-Verlag, Berlin, 1989. \smallskip
\item{[DJS]} David G., Journ\'e J.-L., Semmes S.,{\sl 
Op\'erateurs de Calder\'on-Zygmund, fonctions para-accr\'etives et interpolation}, 
Rev. Mat. Iberoamericana {\bf  1} (1985),  1--56. 
\smallskip
\item{[DvC]} Demuth M., van Casteren J.M., {\sl Stochastic Spectral Theory 
of 
Selfadjoint Feller Operators},
 Birkh\"auser Verlag, Basel, 2000.\smallskip
\item{[FH]}  Fefferman C., Hirachi K., {\sl Ambient metric construction of $Q$-curvature in conformal and CR geometries}, Math. Res. Lett. {\bf10} (2003),  819-831.\smallskip
\item{[Fo1]}  Folland G.B., {\sl  A fundamental solution for a subelliptic operator}, Bull. Amer. Math. Soc. {\bf79} (1973), 373-376.\smallskip
\item{[Fo2]} Folland G.B., {\sl  Subelliptic estimates and function spaces on nilpotent Lie groups}, Ark. Mat. {\bf13} (1975), 161-207.\smallskip
\item{[F]}  Fontana L.,  {\sl Sharp borderline Sobolev inequalities on compact Riemannian manifolds}, Comment. Math. Helv. 
{\bf68} (1993), 415--454.
\smallskip
\item{[FM]}  Fontana L., Morpurgo C., {\sl Adams inequalities in measure spaces},  Adv. Math. {\bf226} (2011), 5066-5119.\smallskip
\item{[FL]} Frank R.L., Lieb  E.H., {\sl Sharp constants in several inequalities on the Heisenberg group}. Ann. of Math., to appear, 	(2010) arXiv:1009.1410v1.\smallskip
\item{[Ge]} Geller D., {\sl The Laplacian and the Kohn Laplacian for the sphere}, J. Differential Geom. {\bf 15} (1980), 417-435.\smallskip 
\item{[GG]} Gover A. R., Graham C.R., {\sl CR invariant powers of the sub-Laplacian}, J. Reine Angew. Math. {\bf 583} (2005), 1-27. \smallskip\item{[Gr]} Graham C.R., {\sl Compatibility operators for degenerated elliptic equations on the ball and Heisenberg group}, Math. Z. {\bf 187} (1984),  289-304.\smallskip
\item {[GJMS]} Graham C.R., Jenne R., Mason L., Sparling G., {\sl 
Conformally invariant powers of the Laplacian, I: existence}, J. London Math. 
Soc. {\bf46} (1992), 557-565.
\smallskip

\item{[H]} Hersch J.,
{\sl  Quatre propri\'et\'es
 isop\'erim\'etrique de
 membranes sph\'erique,  homog\`enes}, 
C. R. Acad. Sci. Paris S\'er. A-B {\bf270} (1970), A1645-A1648. 
\smallskip
\item{[Hi]} Hirachi K., {\sl Scalar pseudo-Hermitian invariants and the Szego kernel on three-dimensional CR manifolds}, Complex geometry (Osaka, 1990), Lecture Notes in Pure and Appl. Math. {\bf 143} (1993), 67-76. \smallskip
\item{[JL1]} Jerison D., Lee J.M., {\sl The Yamabe problem on CR manifolds}, J. Differential Geom. {\bf25} (1987), 167-197.
\smallskip
\item{[JL2]} Jerison D., Lee J.M., {\sl Extremals for the Sobolev inequality on the Heisenberg group and the CR Yamabe problem},  J. Amer. Math. Soc. {\bf1} (1988), 1-13.\smallskip
\item{[JW]}  Johnson K.D., Wallach N.R., {\sl  Composition series and intertwining operators for the spherical principal series. I}, Trans. Amer. Math. Soc. {\bf229} (1977), 137-173.\smallskip
\item{[L]} Lieb E.H., {\sl Sharp constants in the 
Hardy-Littlewood-Sobolev and related inequalities}, 
  Ann. of Math. {\bf118} (1983), 349-374.\smallskip
\item{[KR1]}  Kor\'anyi A., Reimann H.M., {\sl  Quasiconformal mappings on the Heisenberg group}, Invent. Math. {\bf80} (1985), 309-338. 
\smallskip
\item{[KR2]} Kor\'anyi A., Reimann H.M., {\sl Foundations for the theory of quasiconformal mappings on the Heisenberg group}, Adv. Math. {\bf111} (1995), 1-87.\smallskip 
\item{[M1]}Morpurgo C.,{\sl The logarithmic
Hardy-Littlewood-Sobolev
inequality and extremals of zeta functions on~$\,S^n$},
Geom. and Funct. Anal. {\bf6} (1996) 146-171.\smallskip
\item{[M2]}  Morpurgo C., {\sl Sharp inequalities for functional integrals and traces of conformally invariant operators},   
 Duke Math. J. {\bf114} (2002), no. 3, 477--553.
\smallskip
\item{[Mos1]} Moser, J. {\sl A sharp form of an inequality by N. Trudinger}, Indiana Univ. Math. J. {\bf20} (1970/71), 1077-1092.\smallskip
\item{[Mos2]}  Moser J., {\sl On a nonlinear problem in differential geometry}, Dynamical systems (Proc. Sympos., Univ. Bahia, Salvador, 1971), 273-280. Academic Press, New York, 1973. 
\smallskip
\item{[Ok]} Okikiolu, K., {\sl Extremals for Logarithmic HLS inequalities on compact manifolds}, arXiv:math/0603717, (2006).\smallskip
\item{[O]} Onofri, {\sl On the positivity of the
 effective action in a theory
 of random surfaces}, Comm. Math. Phys. {\bf 86} (1982),
 321-326.\smallskip
\item{[OPS]} Osgood B., Phillips R., Sarnak P., {\sl
 Extremals of 
determinants of Laplacians}, J. Funct. Anal. {\bf80} (1988), 
148-211.\smallskip
\item{[Sh]} Showalter R.E., {\sl Hilbert space methods for partial differential equations}, Monographs and Studies in Mathematics, Vol. 1. Pitman, London-San Francisco, Calif.-Melbourne, 1977.\smallskip
\item{[St]} Stanton N.K., {\sl Spectral invariants of CR manifolds}, Michigan Math. J. {\bf 36} (1989),  267-288.\smallskip 
\item{[Ta]} Talenti G., {\sl
Best constant in Sobolev inequality}, 
Ann. Mat. Pura Appl. (0.4) {\bf 110} (1976), 353-372. 
\smallskip

\item{[Tr]} Trudinger N.S., {\sl 
On imbeddings into Orlicz spaces and some applications}, 
J. Math. Mech. {\bf17} (1967) 473--483. 
\smallskip
\item{[VK]}  Vilenkin N.Ja., Klimyk A.U., {\sl Representation of Lie groups and special functions. Vol. 2. Class I representations, special functions, and integral transforms}, Mathematics and its Applications (Soviet Series), 74. Kluwer,  1993.
\smallskip
\item{[Zhu]} Zhu K., {\sl  Spaces of holomorphic functions in the unit ball},  Graduate Texts in Mathematics, 226. Springer-Verlag, New York, 2005.

\bigskip\bigskip\bigskip
\noin Thomas P. Branson (deceased)

\noin Department of Mathematics

\noin University of Iowa

\noin Iowa City, IA 52242

\noin USA
\noin 
\bigskip
\noin Luigi Fontana \hskip19em Carlo Morpurgo

\noin Dipartimento di Matematica e Applicazioni \hskip6em Department of Mathematics 
 
\noin Universit\'a di Milano-Bicocca\hskip 13em University of Missouri, Columbia

\noin Via Cozzi, 53 \hskip 19.3em Columbia, Missouri 65211

\noin 20125 Milano - Italy\hskip 16.6em USA 
\smallskip\noin luigi.fontana@unimib.it\hskip 15.3em morpurgoc@missouri.edu
\end